\documentclass[12pt,reqno]{amsart}
\usepackage{amsfonts,amssymb,latexsym,amsmath,amsthm,color,dsfont}
\usepackage{listings}
\usepackage{textcomp}
\usepackage{bm}
\usepackage{diagbox} 
\usepackage{enumerate,enumitem,cite}
\usepackage{stmaryrd}
\usepackage{fullpage}
\usepackage{multirow}
\usepackage{makecell,array}
\usepackage{tablefootnote}

\usepackage{tabularx}
\usepackage{graphicx, amsmath, amssymb, amsthm,anysize,float,booktabs,geometry}
\usepackage{pdflscape}
\usepackage{caption} 

\usepackage{hyperref}

\hypersetup{citecolor=red, linkcolor=blue, colorlinks=true}

\marginsize{2.54cm}{2.54cm}{2.54cm}{2.54cm}
\allowdisplaybreaks

\newcommand\F{\mathbb{F}}
\newcommand\Z{\mathbb{Z}}
\newcommand{\Q}{\mathbb{Q}}
\newcommand{\N}{\mathbb{N}}

\newcommand{\mc}{\mathcal}

\newcommand{\rank}{\mathrm{rank}}

\newcommand{\wt}{\mathrm{wt}}
\newcommand{\AG}{\mathrm{AG}}
\newcommand{\PG}{\mathrm{PG}}

\newcommand{\supp}{\mathrm{supp}}

\theoremstyle{plain}
\newtheorem{theorem}{Theorem}[section]
\newtheorem{problem}[theorem]{Problem}
\newtheorem{lemma}[theorem]{Lemma}
\newtheorem{corollary}[theorem]{Corollary}

\newtheorem{proposition}[theorem]{Proposition}
\newtheorem{conjecture}[theorem]{Conjecture}

\theoremstyle{definition}
\newtheorem{remark}[theorem]{Remark}
\newtheorem{definition}[theorem]{Definition}
\newtheorem{example}[theorem]{Example}

\numberwithin{equation}{section}

\def\<{\langle}
\def\>{\rangle}

\title{On the maximum size of $(a,b)$-town$\pmod{k}$ families}
\author{Hanlin Zou}

\address{Hanlin Zou, School of Mathematics and Statistics, Yunnan University, Kunming 650091, China}
\email{zouhanlin@ynu.edu.cn}

\keywords{Set system; Oddtown; Eventown.\\
\indent 2020 Mathematics Subject Classification: 05D05}

\begin{document}
\maketitle
\vspace{-0.9cm}

\begin{abstract}
For integers $n \geq k \geq 2$ and $0 \leq a,b \leq k-1$, let
$m_{k,n}(a,b)$ denote the maximum size of an $(a,b)$-town
(mod $k$) family of an $n$-element set, a collection of subsets of whose cardinalities are congruent to $a$ modulo $k$ and whose pairwise intersections are congruent to $b$ modulo $k$. This notion generalizes the classical Oddtown and Eventown problems.

We prove that $m_{k,n}(a,b)\leq n$ whenever
$a\not\equiv b\pmod{k}$, thereby resolving a conjecture of
Veselinov and Marinov. We also disprove another conjecture of theirs by showing that
$m_{3,11}(2,2)>m_{3,11}(1,1)$. For the diagonal case $a\equiv b\pmod{k}$, we establish the general bound
$m_{k,n}(a,a)\leq 2^{\lfloor n/2\rfloor}$ and completely
determine when equality holds. We further obtain improved bounds
and exact values in several special cases. The proofs combine
characteristic-zero linear algebra with methods from coding
theory and finite geometry.
\end{abstract}

\section{Introduction}\label{sec_intro}

The classical Oddtown and Eventown problems are among the most elegant applications of linear algebra in combinatorics. A family $\mc{F}$ of subsets of $[n]:=\{1,\ldots,n\}$ is called an {\it Oddtown} if every member of $\mc{F}$ has odd size and any two distinct members of $\mc{F}$ intersect in an even number of elements; an {\it Eventown} is defined similarly with all members of $\mc{F}$ of even size. It is a classic result that the maximum size of an Oddtown is $n$ (first proved by Berlekamp \cite{Ber69}), while the maximum size of an Eventown is $2^{\lfloor n/2\rfloor}$ (proved independently by Berlekamp \cite{Ber69} and Graver \cite{Gra75}). The proofs rely on simple but powerful linear algebra arguments over finite fields. Over the years, several generalizations of these problems have been studied, see \cite{AC22, BCZ25, FO83, FT16, GST22, JON25, Lin26, ONe23, ONV22, SV18,  WZG25, WZZG25} and the references therein.

Recently, Veselinov and Marinov \cite{VM25} studied the following generalization which they call $(a,b)$-town$\pmod{k}$ families of subsets of $[n]$.

\begin{definition}
Let $a,b,k,n$ be nonnegative integers with $n\geq k\geq 2$. A family $\mc{F}$ of subsets of $[n]$ is called an {\it $(a,b)$-town$\pmod{k}$ family} if for any $A, B\in \mc{F}$ we have
\[
|A\cap B|\equiv \begin{cases}
a\pmod{k},&\text{if }A=B,\\
b\pmod{k},&\text{if }A\neq B.
\end{cases}\]
\end{definition}

Since the intersection conditions depend only on the residues of $a$ and $b$ modulo $k$, we may always assume $a,b\in\{0,1,\ldots,k-1\}$; the only exception is when we construct families where $a$ or $b$ may exceed $k-1$ for convenience.

A natural question is: given parameters $a,b,k,n$, what is the largest possible size of an $(a,b)$-town$\pmod{k}$ family on $[n]$? We denote this maximum by $m_{k,n}(a,b)$. 

When $k=2$, the question is fully understood. 
The values of $m_{2,n}(a,b)$ for all $a,b\in\{0,1\}$ are as follows.
\begin{align*}
m_{2,n}(0,0)=2^{\lfloor n/2\rfloor},& \qquad m_{2,n}(1,0)=n,\\
m_{2,n}(0,1)=\begin{cases}
n, & \text{if $n$ is odd},\\
n-1, & \text{if $n$ is even},
\end{cases}&\qquad
m_{2,n}(1,1)=\begin{cases}
2^{\lfloor (n-1)/2\rfloor}, & \text{if $n$ is odd},\\
2^{\lfloor (n-2)/2\rfloor}, & \text{if $n$ is even}.
\end{cases}
\end{align*}
The $(1,0)$-town and $(0,0)$-town families modulo $2$ are precisely the classical Oddtown and Eventown, respectively. The $(0,1)$ and $(1,1)$ cases follow from \cite[Ex.~1.1.5]{BF22} and \cite[Theorem~A.1]{JON25}, respectively.

When $k\geq 3$, on the other hand, the question is widely open. The first systematic study was undertaken by Veselinov and Marinov \cite{VM25}, who obtained the following bounds for prime modulus $p$.
\begin{equation*}
m_{p,n}(a,b)\leq\begin{cases}
n,&\text{if }a\not\equiv  b\pmod{p},\\
2^{\lfloor (n+1)/2\rfloor}, &\text{if }a\equiv b\pmod{p}.
\end{cases}
\end{equation*}
In the case $a\not\equiv b\pmod{p}$, they obtained a better bound $m_{p,n}(a,b)\leq n-1$ under additional arithmetic conditions on $a,b,n$. Based on these results and further computational evidence, Veselinov and Marinov \cite{VM25} proposed the following two conjectures.

\begin{conjecture}[\cite{VM25}]\label{conj_aa}
Let $a,b,k,n$ be integers with $n\geq k\geq 2$ and $0\leq a,b\leq k-1$. If $a<b$, then $m_{k,n}(a,a)\geq m_{k,n}(b,b)$.
\end{conjecture}

\begin{conjecture}[\cite{VM25}]\label{conj_ab}
Let $a,b,k,n$ be integers with $n\geq k\geq 2$ and $0\leq a,b\leq k-1$. If $a\neq b$, then $m_{k,n}(a,b)\leq n$.
\end{conjecture}

The two conjectures turn out to have different fates. 
We disprove Conjecture \ref{conj_aa} by constructing a $(2,2)$-town$\pmod 3$ family on $[11]$ of size $12$ and proving that no $(1,1)$-town$\pmod{3}$ family on $[11]$ has this size. 
On the other hand, we prove Conjecture \ref{conj_ab} in full generality.
In the diagonal case $a=b$, we also prove a sharp general upper bound and determine all equality cases.
Our main results are as follows.

\begin{theorem}\label{thm_main_aa}
Let $a,k,n$ be integers with $n\geq k\geq 2$ and $0\leq a\leq k-1$. Then
\[m_{k,n}(a,a)\leq 2^{\lfloor\frac{n}{2}\rfloor}.\]
The equality holds if and only if
\begin{enumerate}[label=\upshape{(\roman*)}]
\item $k=2$, $a=0$; or
\item $k=2$, $a=1$, and $n$ is odd; or
\item $n=k=3$ and $a=0$.
\end{enumerate}
Moreover, when $n=k=3$, the only $(0,0)$-town$\pmod{k}$ family of size $2^{\lfloor \frac{n}{2}\rfloor}$ is $\mc{F}=\{\varnothing, \{1,2,3\}\}$.
\end{theorem}

\begin{theorem}\label{thm_main_ab}
Let $a,b,k,n$ be nonnegative integers with $n\geq k\geq 2$ and $0\leq a,b\leq k-1$. If $a\neq b$, then 
\[m_{k,n}(a,b)\leq n.\]
\end{theorem}

In several cases, we also obtain a number of better bounds and exact values beyond the two main theorems. 
In the diagonal case, for $k=3$, we use Type III codes to prove that
$m_{3,n}(a,a)<2.71\cdot 1.3732^n,$
which improves the exponential rate in the general bound $2^{\lfloor n/2\rfloor}$. 
In the off-diagonal case, we determine $m_{k,n}(t,t-1)$ exactly in several congruence classes of $n$. 
For instance,
$m_{k,n}(2,1)=n-1$ if $n\equiv 2\pmod{k}$, 
$m_{3,n}(0,2)=n-2$ if $n\equiv 0\pmod{3}$,
and
$m_{4,n}(0,3)=n-3$ if $n\equiv 0\pmod{4}$.
A more complete list is given in Table \ref{tab_results} at the end of Section \ref{sec_ab}.

We briefly explain the main difficulty in proving Theorem~\ref{thm_main_ab}. When $k$ is a prime, it is often possible to translate congruence conditions modulo $k$ into orthogonality relations over the finite field $\F_k$, and then apply linear algebraic arguments. For composite moduli, however, this approach typically loses essential arithmetic information, since passing from congruences modulo $k$ to a field necessarily discards the interactions among different prime power divisors of $k$.

This phenomenon has already appeared in other extremal set-theoretic problems. A notable example is the $\ell$-Oddtown problem, which asks for the maximum size of a family $\mc{F}$ such that $|A|\not\equiv 0 \pmod{\ell}$ for every $A\in\mc{F}$ and $|A\cap B|\equiv 0 \pmod{\ell}$ for distinct $A,B\in\mc{F}$. While the prime-modulus case is well understood, obtaining strong upper bounds for composite moduli remains notoriously difficult. Only recently did Bukh, Chao, and Zheng \cite{BCZ25} obtain a modest improvement in this direction. These difficulties suggest that new ideas are needed to handle composite moduli.

Our approach avoids reducing the problem to a finite field. Instead, we work directly over characteristic zero and develop a linear independence criterion that preserves congruence information modulo composite numbers. This enables us to establish the sharp bound in Theorem~\ref{thm_main_ab} for all moduli.

The remainder of the paper is organized as follows. 
Section \ref{sec_prelim} contains preliminary material, including the characteristic-zero linear independence criterion, background on self-dual codes, and basic constructions. 
Section \ref{sec_aa} treats the diagonal case $a\equiv b\pmod{k}$, including the counterexample to Conjecture \ref{conj_aa}, the proof of Theorem \ref{thm_main_aa}, and the improved upper bound of $m_{3,n}(a,a)$. 
Section \ref{sec_ab} treats the off-diagonal case $a\not\equiv b\pmod{k}$, including the proof of Theorem \ref{thm_main_ab} and the special exact value results. 
Section \ref{sec_discuss} presents several open problems and future directions.

\section{Preliminaries}\label{sec_prelim}
\subsection{Examples and basic properties of $(a,b)$-town (mod $k$) families}

In this subsection we present several constructions of $(a,b)$-town (mod $k$) 
families, and establish some elementary properties of the maximum size 
$m_{k,n}(a,b)$. 

\begin{example}\label{eg_ab}
Let $a,b,k,n$ be nonnegative integers with $n\geq k\geq 2$ and $a\not\equiv b\pmod{k}$. Suppose that $a-b\equiv s\pmod{k}$ for some $s\in \{1,\ldots,k-1\}$ and let $m=\lfloor\frac{n-b}{s}\rfloor$. For $i\in [m]$, define a subset $F_i$ of $[n]$ by
\[F_i=[b]\cup\{b+(i-1)s+1,\ldots,b+is\}.\]
Then $F_1,\ldots,F_m$ form an $(a,b)$-town$\pmod{k}$ family on $[n]$ with $m=\lfloor\frac{n-b}{s}\rfloor\geq \lfloor\frac{n-b}{k-1}\rfloor$.
\end{example}

\begin{example}\label{eg_aa}
Let $a,k,n$ be nonnegative integers with $n\geq k\geq 2$ and $0\leq a\leq k-1$. Let $m=\lfloor\frac{n-a}{k}\rfloor$. For $i\in [m]$, define a subset $A_i$ of $[n]$ by
\[A_i=\{a+(i-1)k+1,\ldots,a+ik\}.\]
Then $A_1,\ldots,A_m$ are pairwise disjoint and they all have size $k$. Let 
\[\mc{F}=\{[a]\cup S: S\subseteq \{A_1,\ldots,A_m\}\}.\]
Then $\mc{F}$ is an $(a,a)$-town$\pmod{k}$ family on $[n]$ with $|\mc{F}|=2^{\lfloor\frac{n-a}{k}\rfloor}$.
\end{example}

\begin{example}\label{eg_design}
Let $\mc{D}=(\mc{P},\mc{B})$ be a $2$-$(v,k,\lambda)$ symmetric design. Then $|\mc{P}|=|\mc{B}|=v$, and the blocks $\mc{B}$ form a $(k,\lambda)$-town$\pmod{\ell}$ family on $\mc{P}$ for any integer $\ell\geq 2$. When $k\not\equiv \lambda\pmod{\ell}$, this family attains the upper bound $m_{\ell,v}(k,\lambda)\leq v$ in Theorem \ref{thm_main_ab}. Infinitely many such symmetric designs exist; in particular, Corollary \ref{cor_pg2q} will provide an explicit infinite family arising from projective planes.
\end{example}

\begin{example}\label{eg_affine}
Let $n\geq 2$, $q$ a prime power, and $V=\mathbb{F}_q^n$ where $\F_q$ is the finite field of order $q$. Denote by $\mathcal{P}$ the set of vectors in $V$ and by $\mathcal{L}$ the set of cosets of $1$-dimensional subspaces of $V$; these are the points and lines of the affine space $\mathrm{AG}(n,q)$. Then $|\mathcal{P}|=q^n$ and $|\mathcal{L}|=\frac{q^{n-1}(q^n-1)}{q-1}$. For each point $P\in\mathcal{P}$, let $F_P$ be the set of lines passing through $P$. Then $|F_P|=\frac{q^n-1}{q-1}$. Since any two distinct points determine a unique line, we have $|F_P\cap F_Q|=1$ for $P\neq Q$. It follows that $\{F_P: P\in\mathcal{P}\}$ is a $(\frac{q^n-1}{q-1},1)$-town$\pmod{k}$ family on $\mathcal{L}$ for any integer $k\geq 2$.
\end{example}

\begin{example}\label{eg_srg}
Let $G$ be a strongly regular graph with parameters $(v,k,\lambda,\mu)$. For each vertex $v$, let $N_v$ denote its neighbors, so $|N_v|=k$. For two vertices $u$ and $v$, we have $|N_u\cap N_v|=\lambda$ if $u$ and $v$ are adjacent, and $\mu$ otherwise. If $\lambda\equiv \mu\pmod{\ell}$ for some integer $\ell\geq 2$, then $\{N_v: v\in V(G)\}$ forms a $(k,\lambda)$-town$\pmod{\ell}$ family on $V(G)$.
\end{example}

\begin{example}[Frankl-Odlyzko construction via Hadamard matrices]\label{eg_FO}
Let $H$ be a Hadamard matrix of order $4k$. For $i=1,\ldots,4k$, define
\begin{align*}
S_i&=\{j\in [4k]: H_{ij}=1\},\\
T_i&=\{j\in [4k]: H_{ij}=-1\}.
\end{align*}
Then $\{S_1,\ldots,S_{4k},T_1,\ldots,T_{4k}\}$ is a $(0,0)$-town$\pmod{4k}$ family on $[4k]$.
\end{example}

Next, we present some elementary properties for $m_{k,n}(a,b)$. The following lemma generalizes a result of Frankl and Odlyzko \cite[Lemma~3]{FO83}: one can combine several $(a,a)$-town$\pmod{k}$ families on disjoint ground 
sets to obtain a larger $(a,a)$-town$\pmod{k}$ family.

\begin{lemma}\label{lem_combine}
Let $n\geq k\geq 2$ and let $n_1,\ldots,n_s$ be positive integers with 
$n_1+\cdots+n_s=n$. Then
\[
m_{k,n}(a_1+\cdots+a_s,a_1+\cdots+a_s)
   \geq m_{k,n_1}(a_1,a_1)\cdots m_{k,n_s}(a_s,a_s).
\]
\end{lemma}
\begin{proof}
It suffices to prove the case $s=2$ and the general case follows by 
induction on $s$. Let $\{A_1,\ldots,A_m\}$ be an 
$(a_1,a_1)$-town$\pmod{k}$ family on $[n_1]$ with 
$m=m_{k,n_1}(a_1,a_1)$, and let $\{B_1,\ldots,B_\ell\}$ 
be an $(a_2,a_2)$-town$\pmod{k}$ family on 
$\{n_1+1,\ldots,n_1+n_2\}$ with $\ell=m_{k,n_2}(a_2,a_2)$. 
For $i\in[m]$ and $j\in[\ell]$, define
\[
C_{ij}=A_i\cup B_j.
\]
Then $|C_{ij}|=|A_i|+|B_j|\equiv a_1+a_2\pmod{k}$. Moreover, for any 
$(i,j)\neq(u,v)$, we have
\[
|C_{ij}\cap C_{uv}|=|A_i\cap A_u|+|B_j\cap B_v|\equiv a_1+a_2\pmod{k}.
\]
Thus $\{C_{ij}:i\in[m],\,j\in[\ell]\}$ is 
an $(a_1+a_2,a_1+a_2)$-town$\pmod{k}$ family on $[n]$ of size 
$m\ell$, and the inequality follows.
\end{proof}

Combining Example \ref{eg_FO} with Lemma \ref{lem_combine} yields the
following lower bound, originally obtained by Frankl and Odlyzko~\cite{FO83}:
\begin{equation}\label{eq_FObound}
m_{k,n}(0,0)\geq (8k)^{\lfloor\frac{n}{4k}\rfloor},
\end{equation}
for any $n$, under the assumption that a Hadamard matrix of order $4k$ exists.

The following proposition provides several simple relations among the maximum sizes for different parameters. They are obtained by elementary set transformations and will be useful in later sections.

\begin{proposition}\label{prop_basic}
Let $a,b,k,n$ be integers with $n\geq k\geq 2$ and $0\leq a\leq b\leq k-1$. The following statements hold.
\begin{enumerate}[label=\upshape{(\roman*)}]
\item $m_{k,n}(a,b)\leq m_{k,n+1}(a+1,b+1)$.
\item $m_{k,n}(k-1,k-1)\leq m_{k,n+1}(0,0)-1$.
\item $m_{k,n}(a,b)= m_{k,n}(n-a,n-2a+b)$. 
\item $m_{k,n}(a,b)\leq m_{k,n}(2a-2b,a-b)+1$.
\end{enumerate}
\end{proposition}
\begin{proof}
\begin{enumerate}[label=\upshape{(\roman*)}]
\item Let $\mc{F}$ be an $(a,b)$-town$\pmod{k}$ family on $[n]$ with $|\mc{F}|=m_{k,n}(a,b)$. Define a new family $\mc{F}^{(1)}$ by
\[\mc{F}^{(1)}=\{F\cup\{n+1\}: F\in\mc{F}\}.\]
Then $\mc{F}^{(1)}$ is an $(a+1,b+1)$-town$\pmod{k}$ family on $[n+1]$ and $|\mc{F}^{(1)}|=|\mc{F}|$. Therefore $m_{k,n}(a,b)=|\mc{F}^{(1)}|\leq m_{k,n+1}(a+1,b+1)$.
\item Let $\mc{F}$ be a $(k-1,k-1)$-town$\pmod{k}$ family on $[n]$ with $|\mc{F}|=m_{k,n}(k-1,k-1)$. Define a new family $\mc{F}^{(2)}$ by
\[\mc{F}^{(2)}=\{\varnothing\}\cup\{F\cup\{n+1\}: F\in\mc{F}\}.\]
Then $\mc{F}^{(2)}$ is a $(0,0)$-town$\pmod{k}$ family on $[n+1]$ and $|\mc{F}^{(2)}|=|\mc{F}|+1$. As a consequence, we have $m_{k,n}(k-1,k-1)=|\mc{F}|=|\mc{F}^{(2)}|-1\leq m_{k,n+1}(0,0)-1$.
\item Let $\mc{F}$ be an $(a,b)$-town$\pmod{k}$ family on $[n]$ with $|\mc{F}|=m_{k,n}(a,b)$. Define a new family $\mc{F}^{(3)}$ by
\[\mc{F}^{(3)}=\{[n]\setminus F: F\in \mc{F}\}.\]
Then $\mc{F}^{(3)}$ is a $(n-a,n-2a+b)$-town$\pmod{k}$ family on $[n]$. Thus $m_{k,n}(a,b)\leq m_{k,n}(n-a,n-2a+b)$. Similarly, we also have $m_{k,n}(n-a,n-2a+b)\leq m_{k,n}(a,b)$.
\item Let $\mc{F}$ be an $(a,a)$-town$\pmod{k}$ family on $[n]$ with $|\mc{F}|=m_{k,n}(a,a)$. Fix $F_0\in \mc{F}$ and define a new family $\mc{F}^{(4)}$ by
\[\mc{F}^{(4)}=\{F_0\triangle F: F\in \mc{F}\setminus\{F_0\}\},\]
where $F_0\triangle F=(F_0\setminus F)\cup (F\setminus F_0)$ is the symmetric difference of $F_0$ and $F$. Then $\mc{F}^{(4)}$ is a $(2a-2b,a-b)$-town$\pmod{k}$ family on $[n]$. Therefore, we have
$m_{k,n}(a,a)=|\mc{F}^{(4)}|+1\leq m_{k,n}(2a-2b,a-b)+1$.
\end{enumerate}
\end{proof}

\subsection{Linear independence over $\Q$}
Given two vectors $\bm{u}=(u_1,\ldots,u_n)$ and $\bm{v}=(v_1,\ldots,v_n)$ over any field, we use $\bm{u}\cdot \bm{v}$ to denote the standard inner product $\bm{u}\cdot\bm{v}=u_1v_1+\cdots+u_nv_n$.

For a subset $F\subseteq[n]$, its \emph{characteristic vector} is the vector $\bm{v}_F\in\{0,1\}^n$ whose $i$-th coordinate is $1$ if $i\in F$ and $0$ otherwise. The map $F\mapsto\bm{v}_F$ translates set-theoretic conditions into linear-algebraic ones: for any $F,G\subseteq[n]$, we have $\bm{v}_F\cdot\bm{v}_G=|F\cap G|$. 

The following lemma gives a general criterion for the characteristic vectors of certain set families to be linearly independent over $\mathbb{Q}$. It will serve as a key tool in Section \ref{sec_ab}.

\begin{lemma}\label{lem_linindep}
Let $a,b,k,n,m$ be nonnegative integers with $n\geq k\geq 2$ and $m\geq 1$. Let $\bm{v}_1,\ldots,\bm{v}_m$ be $0\text{-}1$ vectors of length $n$ such that 
\[
\bm{v}_i\cdot \bm{v}_j\equiv \begin{cases}
a\pmod{k},&\text{if }i=j,\\
b\pmod{k},&\text{if }i
\neq j.
\end{cases}
\]
If the parameters satisfy the condition
\begin{equation}
\gcd\Big(a(a-b), (a-b)(a-b+mb), k\Big) < k,
\end{equation}
then $\bm{v}_1,\ldots,\bm{v}_m$ are linearly independent over $\Q$.
\end{lemma}

\begin{proof}
Write $d:=\gcd\Big(a(a-b), (a-b)(a-b+mb), k\Big)$. Suppose to the contrary that $\lambda_1\bm{v}_1+\cdots+\lambda_m\bm{v}_m=\bm{0}$ is a primitive relation over $\Z$ with $\gcd(\lambda_1,\ldots,\lambda_m)=1$. 
Taking the inner product of this relation with $\bm{v}_i$ and with the all-one vector $\bm{1}$ respectively, and letting $S = \sum_{j=1}^m \lambda_j$, we obtain
\begin{align}
(a-b)\lambda_i + bS &\equiv 0 \pmod{k}, \quad \forall i \in [m], \label{eq_u1}\\
aS &\equiv 0 \pmod{k}. \label{eq_u2}
\end{align}
Summing \eqref{eq_u1} over all $i \in [m]$ yields 
\[(a-b+mb)S \equiv 0 \pmod k.\] 

Next, multiplying \eqref{eq_u1} by $a$ and using $aS \equiv 0 \pmod k$, we get 
\begin{equation}\label{eq_linindep1}
a(a-b)\lambda_i \equiv 0 \pmod k.
\end{equation}
Similarly, multiplying \eqref{eq_u1} by $(a-b+mb)$ and using $(a-b+mb)S \equiv 0 \pmod k$, we get 
\begin{equation}\label{eq_linindep2}
(a-b)(a-b+mb)\lambda_i \equiv 0 \pmod k.
\end{equation}
Now, \eqref{eq_linindep1} and \eqref{eq_linindep2} imply that $k' \mid \lambda_i$ for all $i$, where $k'=\frac{k}{d}$. Since $d<k$, then $k'>1$. This contradicts our assumption that $\gcd(\lambda_1,\ldots,\lambda_m)=1$. Therefore $\bm{v}_1,\ldots,\bm{v}_m$ are linearly independent over $\Q$.
\end{proof}

We will also need the following elementary bound on the number of $\{0,1\}$-vectors in a subspace, due to Odlyzko.

\begin{lemma}[{\cite[Theorem 2]{Odl81}}]\label{lem_Odl81}
Let $\F$ be a field and $V\subset \F^n$ be a subspace of dimension $d$. Then $|V\cap \{0,1\}^n|\leq 2^d$.
\end{lemma}

\subsection{Self-dual codes}

Let $q$ be a prime power and $\mathbb{F}_q$ the finite field of order $q$. A {\it linear code} is a subspace $\mathcal{C}$ of the vector space $\mathbb{F}_q^n$. The {\it dual code} of $\mathcal{C}$ is the orthogonal complement $\mathcal{C}^{\perp}=\{\bm{x}\in\mathbb{F}_q^n : \bm{x}\cdot \bm{y} = 0 \text{ for all } \bm{y}\in\mathcal{C}\}$, 
where $\bm{x}\cdot \bm{y}=\sum_{i=1}^n x_i y_i$ is the standard bilinear form. The vectors in $\mathcal{C}$ are called {\it codewords}, and the {\it weight} $\operatorname{wt}(\boldsymbol{x})$ of a codeword $\boldsymbol{x}\in\mathcal{C}$ is the number of nonzero coordinates of $\boldsymbol{x}$. Codes over $\mathbb{F}_2$ and $\F_3$ are called {\it binary codes} and {\it ternary codes}, respectively. 

A linear code $\mathcal{C}$ is {\it self-dual} if $\mathcal{C} = \mathcal{C}^{\perp}$. A binary self-dual code in which all weights are divisible by $4$ is called a {\it Type II} code. A ternary self-dual code in which all weights are divisible by $3$ is called a {\it Type III} code. Type II and Type III codes are rare in the sense that they exist only for lengths satisfying certain congruences, as shown in the following result.

\begin{lemma}[{\cite[Corollary 9.2.2]{HP03}}]\label{lem_type2}
Type II codes of length $n$ exist if and only if $n \equiv 0 \pmod{8}$; and Type III codes exist if and only if $n\equiv 0\pmod{4}$. 
\end{lemma}

Given a linear code $\mc{C}$ of length $n$, the (two variable) {\it weight enumerator} is 
\[W_{\mc{C}}(x,y)=\sum_{i}A_i(\mc{C})x^iy^{n-i},\]
where $A_i(\mc{C})$ is the number of codewords in $\mc{C}$ with weight $i$. For Type III codes, the weight enumerators are described by Gleason's theorem which we state below.

\begin{lemma}[{\cite[Theorem 9.2.1]{HP03}}]\label{lem_Gle}
Let $\Phi_4(x,y)=y^4+8x^3y$ and $\Phi_{12}(x,y)=y^{12}+264x^6y^6+440x^9y^3+24x^{12}$. Then the weight enumerator of a Type III code $\mc{C}$ is
\[W_{\mc{C}}(x,y)=\sum_{i=0}^{\lfloor\frac{n}{12}\rfloor}a_i\Phi_4(x,y)^{\frac{n}{4}-3i}\Phi_{12}(x,y)^i,\]
where all $a_i$'s are rational and $\sum_i a_i=1$. 
\end{lemma}

This provides an estimate of $W_{\mc{C}}(1,2)$, which in turn yields a tighter upper bound for $m_{3,n}(0,0)$.

\begin{corollary}\label{cor_wt12}
Let $\mc{C}$ be a Type III code of length $n$. Then
\[W_{\mc{C}}(1,2)<34^{\frac{n}{4}}.\]
\end{corollary}
\begin{proof}
We deduce from Lemma \ref{lem_Gle} that
\begin{align*}
W_{\mc{C}}(1,2)&=\sum_{i=0}^{\lfloor\frac{n}{12}\rfloor}a_i\Phi_4(1,2)^{\frac{n}{4}-3i}\Phi_{12}(1,2)^i=\sum_{i=0}^{\lfloor\frac{n}{12}\rfloor}a_i34^{\frac{n}{4}-3i}24536^i\\
&<\sum_{i=0}^{\lfloor\frac{n}{12}\rfloor}a_i34^{\frac{n}{4}}=34^{\frac{n}{4}}\sum_{i=0}^{\lfloor\frac{n}{12}\rfloor}a_i=34^{\frac{n}{4}}.
\end{align*}
The last equality holds since $\sum_i a_i=1$. 
\end{proof}

We will also need a basic character sum identity over $\F_3$ for improving the bound for $m_{3,n}(0,0)$ in Subsection \ref{sub_m3naa}. 
Let $\mc{C}$ be a ternary code of length $n$. Define the function $\mathds{1}_{\mc{C}}$ on $\F_3^n$ by
\begin{equation}\label{eq_indi}
\mathds{1}_{\mc{C}}(\bm{x}) = \frac{1}{|\mc{C}^\perp|}\sum_{\bm{u}\in \mc{C}^\perp}\zeta_3^{\bm{x}\cdot\bm{u}}, \qquad \bm{x}\in\F_3^n,
\end{equation}
where $\zeta_3$ is a primitive third root of unity. We claim that $\mathds{1}_{\mc{C}}$ is the indicator function of $\mc{C}$, i.e.,
\begin{equation}\label{eq_char}
\mathds{1}_{\mc{C}}(\bm{x}) = 
\begin{cases}
1, & \text{if } \bm{x}\in\mc{C},\\
0, & \text{otherwise}.
\end{cases}
\end{equation}
To verify this, set $S(\bm{x}) = \sum_{\bm{u}\in \mc{C}^\perp}\zeta_3^{\bm{x}\cdot\bm{u}}$. 
If $\bm{x}\in\mc{C}$, then $\bm{x}\cdot\bm{u}=0$ for every $\bm{u}\in\mc{C}^\perp$, so $S(\bm{x}) = |\mc{C}^\perp|$.
If $\bm{x}\notin\mc{C}$, then since $\mathcal{C} = (\mc{C}^\perp)^\perp$, there exists $\bm{u}_0\in\mc{C}^\perp$ such that $\bm{x}\cdot\bm{u}_0 \neq 0$. 
Set $\alpha = \zeta_3^{\bm{x}\cdot\bm{u}_0} \neq 1$. 
Then
\[
S(\bm{x}) = \sum_{\bm{u}\in \mathcal{C}^\perp}\zeta_3^{\bm{x}\cdot\bm{u}} 
= \sum_{\bm{u}\in \mathcal{C}^\perp}\zeta_3^{\bm{x}\cdot(\bm{u}+\bm{u}_0)} 
= \alpha\, S(\bm{x}),
\]
which forces $S(\bm{x})=0$ because $\alpha\neq 1$. 
Thus $\mathds{1}_{\mathcal{C}}(\bm{x}) = \frac{S(\bm{x})}{|\mathcal{C}^\perp|}$ satisfies \eqref{eq_char}.

\begin{remark}\label{rmk_tispace}
We recall some terminology from finite geometry that will be used in our study of $(0,0)$-town families in Subsection \ref{sub_mknaa}. Let $V=\mathbb{F}_p^n$ be equipped with the standard bilinear form $\bm{x}\cdot \bm{y}=\sum_{i=1}^n x_i y_i$, which is nondegenerate for all primes $p$.
A subspace $W\subseteq V$ is called \emph{totally isotropic} if $\bm{x}\cdot \bm{y}=0$ for all $\bm{x},\bm{y}\in W$.
In this language, a self-orthogonal code is a totally isotropic subspace, and a self-dual code is a maximal totally isotropic subspace.
Witt's theorem states that a maximal totally isotropic subspace of $V$ has dimension $\lfloor \frac{n}{2}\rfloor$. Consequently, a self-dual code of length $n$ can only exist when $n$ is even.
\end{remark}

\subsection{Sums of binomial coefficients over residue classes}

To establish the upper bounds for $m_{k,n}(a,b)$ for certain parameters, it is essential to explicitly evaluate the number of vectors in $\F_2^n$ whose weights satisfy certain modular congruences. This subsection provides the necessary technical lemmas regarding the sums of binomial coefficients over residue classes modulo 4, which are central to our counting arguments in Subsection \ref{sub_betterab}.

For $n,k\in \N^*$ and $0\leq r\leq k-1$, we define the partial binomial sum
\[T(n,r,k):=\sum_{\substack{j=0\\j\equiv r\!\!\!\!\pmod{k}}}^n\binom{n}{j}.\] 
A classical result, often attributed to Ramus \cite{Ram1834}, provides a closed-form expression for $T(n,r,k)$ using roots of unity, that is,
\begin{equation}\label{eq_Ramus}
T(n,r,k)=\frac{1}{k}\sum_{j=0}^{k-1}\omega_k^{(k-j)r}(1+\omega_k^j)^n,
\end{equation}
where $\omega_k=e^{\frac{2\pi i}{k}}$ is a primitive $k$-th root of unity. For our purpose, we focus on the case $k=4$. The following lemma provides an explicit evaluation.

\begin{lemma}\label{lem_Tmr4}
For $n\in \N^*$ and $0\leq r\leq 3$, the sum $T(n,r,4)$ is given by
\[T(n,r,4)=2^{n-2}+2^{\frac{n-2}{2}}\cos\left(\frac{(n-2r)\pi}{4}\right).\]
\end{lemma}
\begin{proof}
Setting $k=4$ and $\omega_4=e^{\frac{\pi i}{2}}=i$ in \eqref{eq_Ramus}, we have
\begin{equation}\label{eq_Tmr4_1}
T(n,r,4)=\frac{1}{4}(2^n+e^{\frac{-r \pi i}{2}}(1+i)^n+e^{\frac{r \pi i}{2}}(1-i)^n).
\end{equation}
Using $1\pm i=\sqrt{2}e^{\frac{\pm\pi i}{4}}$, we have
\begin{equation}\label{eq_Tmr4_2}
e^{\frac{-r \pi i}{2}}(1+i)^n+e^{\frac{r \pi i}{2}}(1-i)^n=2^{\frac{n}{2}}(e^{\frac{(n-2r)\pi i}{4}}+e^{-\frac{(n-2r)\pi i}{4}})=2^{\frac{n+2}{2}}\cos\left(\frac{(n-2r)\pi}{4}\right).
\end{equation}
The desired form of $T(n,r,4)$ then follows by combining \eqref{eq_Tmr4_1} and \eqref{eq_Tmr4_2}.
\end{proof}

\section{The case $a\equiv b\pmod{k}$}\label{sec_aa}

\subsection{A counterexample to Conjecture \ref{conj_aa}}\label{subsec_conj_aa}

In this subsection, we present a counterexample to Conjecture \ref{conj_aa}. More precisely, we prove the following result.
\begin{proposition}
We have $m_{3,11}(1,1)<m_{3,11}(2,2)$.
\end{proposition}
\begin{proof}
First of all, it is easy to verify that 
\begin{align*}
&\{\{ 2, 3, 4, 7, 9 \}, \{ 3, 7, 8, 10, 11 \}, \{ 4, 5, 8, 9, 10 \}, \{ 1, 3, 6, 8, 9 \},\\
&\{ 5, 6, 7, 9, 11 \}, \{ 2, 4, 6, 8, 11 \}, \{ 2, 3, 5, 6, 10 \}, \{ 1, 4, 6, 7, 10 \}, \\
&\{ 1, 2, 5, 7, 8 \}, \{ 1, 2, 9, 10, 11 \}, \{1, 3, 4, 5, 11 \}, \{ 1, 2, 3, 4, 5, 6, 7, 8, 9, 10, 11 \} \}
\end{align*}
is a $(2,2)$-town$\pmod{3}$ family on $[11]$. Therefore 
\[m_{3,11}(2,2)\geq 12.\]

In the rest of this proof, we show that $m_{3,11}(1,1)<12$. Let $\mc{F}$ be a $(1,1)$-town$\pmod{3}$ family on $[11]$. 
Let $\mc{A}=\{A\in\mc{F}:|A|=4\}$ and $\mc{B}=\{B\in\mc{F}:|B|=7\}$, and set $a=|\mc{A}|$, $b=|\mc{B}|$. For each $i\in[11]$, define
\[
x_i = |\{A\in\mc{A}: i\in A\}|.
\]
Counting the pairs $(i,A)$ with $i\in A$ and $A\in\mathcal{A}$, and the pairs $(i,\{A_1,A_2\})$ with $i\in A_1\cap A_2$ for distinct $A_1,A_2\in\mathcal{A}$, we obtain
\begin{equation}\label{eq:sumx}
\sum_{i=1}^{11} x_i = 4a,\quad \text{ and } \sum_{i=1}^{11} \binom{x_i}{2} = \binom{a}{2},
\end{equation}
respectively. Combining these with $\sum_i x_i^2 = 2\sum_i \binom{x_i}{2} + \sum_i x_i$ and the Cauchy--Schwarz inequality $(\sum_i x_i)^2 \leq 11\sum_i x_i^2$ yields $16a^2 \leq 11(a^2+3a)$, hence 
\[a \leq \lfloor 33/5\rfloor = 6.\]
A similar argument gives $b\leq 6$. 

{\bf Case 1.} Assume that $|F|\in\{4,7\}$ for each $F\in \mc{F}$.

In this case, we have $|\mc{F}| = a+b \leq 12$, with equality only possible if $a=b=6$.
Suppose for contradiction that $a=b=6$. Let $t$ be the number of pairs $(A,B)$ with $A\in\mathcal{A}$, $B\in\mathcal{B}$ and $A\subseteq B$. Fix $A\in\mathcal{A}$. Suppose that $A$ is contained in $r(A)$ members of $\mathcal{B}$, say $B_1,\ldots,B_{r(A)}$. Then $B_i\cap B_j=A$ for $i\neq j$. Consequently, we have
\[3r(A)=\sum_{i=1}^{r(A)}|B_i\setminus A|\leq |[11]\setminus A|=7.\]
Thus $r(A)\leq 2$ and it follows that 
\begin{equation}\label{eq:t-ub}
t = \sum_{A\in\mc{A}}r(A) \leq 2a=12.
\end{equation}

Now let $z_i$ be the number of members of $\mathcal{F}=\mathcal{A}\cup\mathcal{B}$ containing $i$. By the same counting arguments as for $x_i$, we have
\[
\sum_i z_i = 4a+7b = 66,\qquad 
\sum_i \binom{z_i}{2} = \binom{a}{2} + 4\binom{b}{2} + 4t + (ab-t)=111+3t.
\]
Combining these with $\sum_{i=1}^{11} z_i^2 = 2\sum_{i=1}^{11} \binom{z_i}{2} + \sum_{i=1}^{11} z_i$ and $(\sum_{i=1}^{11} z_i)^2 \leq 11\sum_{i=1}^{11} z_i^2$, we obtain $66^2 \leq 11(288+6t)$. So $t \geq 18$. This contradicts $t\leq 12$ from \eqref{eq:t-ub}. Hence $a=b=6$ is impossible, and we must have $|\mc{F}|<12$.

{\bf Case 2.} Assume that $\mc{F}$ contains $1$-subsets or $10$-subsets.

Note that $\mc{F}$ can contain at most one $1$-subset, and similarly at most one $10$-subset. By taking complements if necessary, we may assume that $\mc{F}$ contains a $1$-subset, say $\{1\}$. Then every $4$-subset in $\mc{F}$ must contain $1$, and any two such $4$-subsets are disjoint outside $\{1\}$. It follows that $a\leq 3$, and consequently, we have $|\mc{F}| \leq 1 + 3 + 6 + 1 = 11$. 

Combining the two cases, we have $m_{3,11}(1,1)\leq 11<12\leq m_{3,11}(2,2)$. 
\end{proof}

\begin{remark}
An exhaustive computer search using Magma~\cite{Magma} confirms the exact values:
$m_{3,11}(1,1)=8$ and $m_{3,11}(2,2)=12$.
\end{remark}

\subsection{General upper bound for $m_{k,n}(a,a)$}\label{sub_mknaa}

This subsection is devoted to the proof of Theorem \ref{thm_main_aa}. We begin with a simple but useful observation.

\begin{lemma}\label{lem:aa-le-00}
For any integers $n\geq k\geq2$, and $0\leq a\leq k-1$, we have
\[
m_{k,n}(a,a) \leq m_{k,n}(0,0).
\]
\end{lemma}
\begin{proof}
Let $\mathcal{F}$ be an $(a,a)$-town$\pmod{k}$ attaining the maximum size. Fix $S\in\mathcal{F}$ and define $\mathcal{F}'=\{S\triangle F:F\in\mathcal{F}\}$. Then $|\mathcal{F}'|=|\mathcal{F}|$. Using the identity $|U\triangle V|=|U|+|V|-2|U\cap V|$ and the fact that $\mathcal{F}$ is an $(a,a)$-town, one readily checks that each member of $\mathcal{F}'$ has size divisible by $k$ and any two distinct members intersect in size divisible by $k$. Thus $\mathcal{F}'$ is a $(0,0)$-town$\pmod{k}$, giving $m_{k,n}(a,a)\leq m_{k,n}(0,0)$.
\end{proof}

\begin{corollary}\label{cor_maa_ub}
For any integers $n\geq k\geq2$ and $0\leq a\leq k-1$, we have
\[
m_{k,n}(a,a)\leq 2^{\lfloor n/2\rfloor}.
\]
\end{corollary}
\begin{proof}
By Lemma~\ref{lem:aa-le-00}, it suffices to show $m_{k,n}(0,0)\leq 2^{\lfloor n/2\rfloor}$. Any $(0,0)$-town$\pmod{k}$ is also a $(0,0)$-town$\pmod{p}$ for any prime divisor $p$ of $k$. Hence the characteristic vectors of a $(0,0)$-town family span a totally isotropic subspace of $\F_p^n$ of dimension $d\leq\lfloor n/2\rfloor$ (see Remark~\ref{rmk_tispace}). By Lemma~\ref{lem_Odl81}, such a subspace contains at most $2^d\leq 2^{\lfloor n/2\rfloor}$ vectors of $\{0,1\}^n$.
\end{proof}

Our next task is to determine when equality can hold in Corollary \ref{cor_maa_ub}.

\begin{proposition}\label{pp_00towm_1}
For integers $n\geq k\geq 2$, if there exists a $(0,0)$-town$\pmod{k}$ family $\mc{F}$ on $[n]$ with $|\mathcal{F}| = 2^{\lfloor n/2 \rfloor}$, then 
\begin{enumerate}[label=\upshape{(\roman*)}]
\item $k$ is a power of 2; or
\item $n=k=3$ and $\mc{F}=\{\varnothing, \{1,2,3\}\}$.
\end{enumerate}
\end{proposition}

\begin{proof}
Suppose that $k$ is not a power of $2$. Then $k$ has an odd prime divisor $p\geq 3$. For $A \in \mathcal{F}$, let $\bm{v}_A \in \{0, 1\}^n$ denote its characteristic vector, regarded as a vector in $\F_p^n$. Then $\bm{v}_A \cdot \bm{v}_B =0$ for all $A, B \in \mathcal{F}$.

Let $V =\langle \{\bm{v}_A : A \in \mathcal{F}\}\rangle_{\F_p}$ be the subspace spanned by these vectors. Then $V$ is totally isotropic, and Remark~\ref{rmk_tispace} gives $\dim(V) \leq \lfloor n/2\rfloor$. By Lemma~\ref{lem_Odl81}, we have
\[
2^{\lfloor n/2 \rfloor}=  |\mathcal{F}|\leq |V\cap \{0,1\}^n|\leq 2^{\dim V}\leq 2^{\lfloor n/2 \rfloor}.
\]
Therefore $\dim(V) = \lfloor n/2 \rfloor$.

Write $d=\lfloor n/2 \rfloor$. Since $\dim(V) = d$, we can choose $d$ coordinate positions $I = \{i_1, \dots, i_d\}$ such that the projection $\pi: V \to \F_p^d$ onto these coordinates is an isomorphism. Consider $S:=\{\bm{v}_A: A\in \mathcal{F}\}$. Since $S \subseteq \{0,1\}^n$, we have $\pi(S) \subseteq \{0,1\}^d$. Because $\pi$ is bijective, we have $|\pi(S)| = |S| = |\mathcal{F}| = 2^d$. The only subset of $\{0,1\}^d$ of size $2^d$ is $\{0,1\}^d$ itself, hence
\begin{equation}\label{eq_00ub_1}
\pi(S) = \{0,1\}^d.
\end{equation}

Let $\bm{e}_1, \ldots, \bm{e}_d$ be the standard basis vectors in $\F_p^d$. Choose $\bm{v}_1, \ldots, \bm{v}_d \in S$ with $\pi(\bm{v}_j) = \bm{e}_j$. 
For distinct $i, j$, the vector $\bm{e}_i + \bm{e}_j$ belongs to $\{0,1\}^d$, so by \eqref{eq_00ub_1} there exists $\bm{w} \in S$ with $\pi(\bm{w}) = \bm{e}_i + \bm{e}_j$. On the other hand, $\pi(\bm{v}_i+\bm{v}_j)=\bm{e}_i+\bm{e}_j$. Since $\pi$ is injective on $V$, we must have $\bm{w} = \bm{v}_i+\bm{v}_j$ in $\F_p^n$.

Now view the equality $\bm{w} = \bm{v}_i+\bm{v}_j$ coordinatewise. If there were a coordinate $k$ with $(\bm{v}_i)_k = (\bm{v}_j)_k = 1$, then $\bm{w}_k = 2$. But $\bm{w}\in S$ implies that $\bm{w}$ is a $\{0,1\}$-vector, a contradiction. Hence $\bm{v}_i$ and $\bm{v}_j$ have disjoint supports. It follows that 
\[\sum_{i=1}^d |\supp(\bm{v}_i)|\leq n.\] 
On the other hand, each $\bm{v}_i$ comes from a member of the $(0,0)$-town family, so $|\supp(\bm{v}_i)| \equiv 0 \pmod p$. Since $\bm{v}_i\neq \bm{0}$, we have $|\supp(\bm{v}_i)| \geq p$. Thus
\[
dp\leq \sum_{i=1}^d |\supp(\bm{v}_i)|\leq n.
\]
This yields $\lfloor n/2\rfloor = d \leq n/p \leq n/3$, which forces $n=3$.

When $n=3$, we have $|\mathcal F|=2$. Since $k$ is not a power of $2$, necessarily $k=3$, in which case the unique $(0,0)$-town family of size $2$ on $[3]$ is $\mc{F}=\{\varnothing, \{1,2,3\}\}$.
\end{proof}

Next, we show that if $k$ is a multiple of $4$, then $m_{k,n}(0,0)<2^{\lfloor n/2\rfloor}$. It suffices to show $m_{4,n}(0,0)<2^{\lfloor n/2\rfloor}$ since any $(0,0)$-town$\pmod{k}$ family is also a $(0,0)$-town$\pmod{4}$ family. 
We treat the two cases according to the parity of $n$ separately in Theorems~\ref{thm_m4n00_even} and~\ref{thm_m4n00_odd}.

\medskip
\noindent\textbf{Notation.} 
For the remainder of this subsection, whenever $\mathcal{F}=\{F_1,\ldots,F_m\}$ is an $(a,b)$-town$\pmod{k}$ family, 
we denote its characteristic vectors by $\bm{v}_1,\ldots,\bm{v}_m$. 
The underlying field will be specified in each context.

\begin{lemma}\label{lem_00mod4_n8}
Suppose that $n$ is even and $k\equiv 0\pmod{4}$. If there exists a $(0,0)$-town$\pmod{k}$ family on $n$ of size $2^{\frac{n}{2}}$, then $n\equiv 0\pmod{8}$.
\end{lemma}
\begin{proof}
Let $\mathcal{F}$ be a $(0,0)$-town$\pmod{k}$ family of size $2^{\frac{n}{2}}$. Let $\mathcal{C}$ be the binary linear code generated by the characteristic vectors of $\mathcal{F}$. Then $\mc{C}$ is self-orthogonal with $|\mc{C}|\geq|\mc{F}|=2^{\frac{n}{2}}$. By Remark \ref{rmk_tispace}, we see that $|\mc{C}|\leq 2^{\frac{n}{2}}$. Thus $\dim(\mc{C})=\frac{n}{2}$ and $\mc{C}=\{\bm{v}_F: F\in\mc{F}\}$. It follows that $\dim(\mc{C}^\perp)=n-\dim(\mc{C})=\frac{n}{2}$ and so $\mc{C}=\mc{C}^\perp$. Therefore $\mc{C}$ is a Type II code, and Lemma \ref{lem_type2} then forces $n\equiv0\pmod{8}$.
\end{proof}

\begin{lemma}\label{lem_m4n30n}
Let $n\geq 4$. Then $m_{4,n}(3,0)<n$.
\end{lemma}
\begin{proof}
Suppose to the contrary that $\mathcal{F}=\{F_1,\ldots,F_n\}$ is a $(3,0)$-town$\pmod{4}$ family on $[n]$. We regard their characteristic vectors $\bm{v}_1,\ldots,\bm{v}_n$ as vectors over $\Q$. Then for all $i,j\in[n]$, we have
\[
\bm{v}_i\cdot \bm{v}_j \equiv
\begin{cases}
3 \pmod{4}, & \text{if }i=j,\\
0 \pmod{4}, & \text{if }i\neq j.
\end{cases}
\]
By Lemma~\ref{lem_linindep} with $(a,b,k,m)=(3,0,4,n)$, we have
\[
\gcd\big(a(a-b),\, (a-b)(a-b+mb),\, k\big) = \gcd(9,\,9,\,4) = 1 < 4,
\]
so $\bm{v}_1,\ldots,\bm{v}_n$ are linearly independent over $\mathbb{Q}$. Hence they form a basis of $\mathbb{Q}^n$. 
Consequently, the all-one vector $\bm{1}$ can be expressed as a rational linear combination of the $\bm{v}_i$'s. 
Clearing denominators and dividing by the greatest common divisor yields integers $c_0,c_1,\ldots,c_n$ with $\gcd(c_0,c_1,\ldots,c_n)=1$ such that
\begin{equation}\label{eq_m4n30_2}
c_1\bm{v}_1+\cdots+c_n\bm{v}_n = c_0\bm{1}.
\end{equation}

For each $j\in [n]$, taking inner products of both sides of \eqref{eq_m4n30_2} with $\bm{v}_j$, we have 
\[3c_j\equiv 3c_0\pmod{4}.\]
Thus 
\[c_j\equiv c_0\pmod{4}.\]

Let $M$ be the $n\times n$ integer matrix whose $i$-th row is $\bm{v}_i$. Then 
\[MM^\top\equiv 3I \pmod{4}.\]
Thus $\det(M)^2\equiv 3^n\equiv (-1)^n\pmod{4}$. If $n$ is odd, then $\det(M)^2\equiv -1\pmod{4}$. This is a contradiction since $-1$ is not a square in $\Z_4$. If $n$ is even, then $\det(M)^2\equiv 1\pmod{4}$, which implies that $\det(M)\equiv 1\text{ or }3\pmod{4}$. Since $1$ and $3$ are both units of $\Z_4$, we see that $M$ is invertible and so
\[M^\top M\equiv 3I \pmod{4}.\]
It follows that 
\begin{equation}\label{eq_m4n30_3}
(\bm{v}_1+\cdots+\bm{v}_n)_j=(\bm{v}_1)_j^2+\cdots+(\bm{v}_n)_j^2=(M^\top M)_{jj}\equiv 3\pmod{4},
\end{equation}
for each $j\in [n]$. Since $c_i\equiv c_0\pmod{4}$ for each $i\in[n]$, we deduce from \eqref{eq_m4n30_2} and \eqref{eq_m4n30_3} that 
\[ c_0\equiv c_0(\bm{v}_1+\cdots+\bm{v}_n)_j\equiv 3c_0\pmod{4},\]
for each $j\in [n]$. This implies that $c_0$ is even and so all the $c_i$'s are even, contradicting the assumption that $\gcd(c_0,\ldots,c_n)=1$. Therefore $m_{4,n}(3,0)<n$. 
\end{proof}

\begin{theorem}\label{thm_m4n00_even}
If $n\geq 4$ is even, then $m_{4,n}(0,0)<2^{\frac{n}{2}}$.
\end{theorem}
\begin{proof}
Write $d=\frac{n}{2}$ and $m=2^{\frac{n}{2}}$. Suppose to the contrary that $\mc{F}=\{F_1,\ldots,F_m\}$ is a $(0,0)$-town$\pmod{4}$ family. We regard their characteristic vectors $\bm{v}_1,\ldots,\bm{v}_m$ as vectors over $\F_2$.  By Lemma \ref{lem_00mod4_n8}, we have $n\equiv 0\pmod{8}$ and $V=\<\bm{v}_1,\ldots,\bm{v}_n\>$ is a maximal totally isotropic subspace of $\F_2^n$ of dimension $d=\frac{n}{2}$. Thus there exist $d$ positions $i_1,\ldots,i_{d}$ such that the projection from $V$ to these $d$ positions is an isomorphism between $V$ and $\F_2^{d}$. Let $\bm{e}_1,\ldots,\bm{e}_d$ be the standard basis for $\F_2^d$. Then there exist $\bm{v}_{k_1},\ldots,\bm{v}_{k_d}$ such that $\pi(\bm{v}_{k_i})=\bm{e}_i$ for $i\in [d]$. So we can write $\bm{v}_{k_i}=(\bm{e}_i, \bm{u}_i)$ where $\bm{u}_i\in\F_2^{n-d}$. Then for $i,j\in [d]$, we have
\[\bm{u}_i\cdot \bm{u}_j\equiv \begin{cases}
3\pmod{4},&\text{if }i=j,\\
0\pmod{4},&\text{if }i\neq j.
\end{cases}\]
Let $A_1,\ldots,A_d$ be the subsets of $[n]\setminus\{i_1,\ldots,i_d\}$ corresponding to $\bm{u}_1,\ldots,\bm{u}_d$, respectively. Then $A_1,\ldots,A_d$ form a $(3,0)$-town$\pmod{4}$  family on $[n]\setminus\{i_1,\ldots,i_d\}$. However, this contradicts Lemma \ref{lem_m4n30n}. Therefore $m_{4,n}(0,0)<2^{\frac{n}{2}}$.
\end{proof}

We now turn to the case where $n$ is odd. Let $n\geq 3$ be odd. Set $d=\frac{n-1}{2}$ and $m=2^d$. Suppose that $\mc{F}=\{F_1,\ldots,F_m\}$ is a $(0,0)$-town$\pmod{4}$ family on $[n]$. We regard their characteristic vectors $\bm{v}_1,\ldots,\bm{v}_m$ as vectors over $\F_2$, and define a binary code by $\mc{C}=\<\bm{v}_1,\ldots,\bm{v}_m\>$.

\begin{lemma}\label{lem_structure_C}
The code $\mc{C}$ has the following properties.
\begin{enumerate}[label=\upshape{(\roman*)}]
\item $\dim(\mc{C})=d$ and $\mc{C}=\{\bm{v}_1,\ldots,\bm{v}_m\}$. Consequently, every codeword of $\mc{C}$ has weight divisible by $4$.
\item $\mc{C}^\perp=\<\bm{1}\>\oplus\mc{C}$.  
Moreover, a vector $\bm{x}\in\mc{C}^\perp$ belongs to $\mc{C}$ if and only if $\wt(\bm{x})\equiv0\pmod{4}$.
\end{enumerate}
\end{lemma}

\begin{proof}
(i) Since $\mc{F}$ is a $(0,0)$-town$\pmod{4}$ family, we see that $\mc{C}$ is self-orthogonal. By Remark \ref{rmk_tispace}, we have $2^d=|\mc{F}|\leq |\mc{C}|\leq 2^d$. Thus $\dim(C)=d$ and $\mc{C}=\{\bm{v}_i: 1\leq i\leq m\}$. It follows that $\wt(\bm{u})\equiv 0\pmod{4}$ for any $\bm{u}\in\mc{C}$. 

(ii) For any $F\in\mc{F}$, we have $\bm{1}\cdot\bm{v}_F = |F|\equiv 0\pmod{4}$, hence $\bm{1}\in\mc{C}^\perp$. Since $\wt(\bm{1})=n\equiv 1\pmod{2}$, we have $\bm{1}\notin\mc{C}$. 
Now $\dim(\mc{C}^\perp) = n-d = d+1$, and $\mc{C}\subsetneq\langle\bm{1}\rangle+\mc{C}\subseteq\mc{C}^\perp$. 
Comparing dimensions forces $\mc{C}^\perp=\langle\bm{1}\rangle\oplus\mc{C}$. 
The last assertion follows because every vector in $\mc{C}^\perp$ has the form $\bm{u}$ or $\bm{1}+\bm{u}$ with $\bm{u}\in\mc{C}$, 
and $\wt(\bm{1}+\bm{u})=n-\wt(\bm{u})\equiv n\not\equiv0\pmod{4}$, while $\wt(\bm{u})\equiv0\pmod{4}$.
\end{proof}

\begin{lemma}\label{lem_Hadamard_closed}
For any $\bm{u},\bm{v}\in\mathcal{C}$, 
their Hadamard product $\bm{u}\circ\bm{v}=(u_1 v_1,\ldots,u_n v_n)$ also lies in $\mathcal{C}$.
\end{lemma}

\begin{proof}
Let $\bm{u},\bm{v},\bm{w}\in\mc{C}$. By Lemma \ref{lem_structure_C}, 
$\bm{u},\bm{v},\bm{w},\bm{u}+\bm{w}$ correspond to members of $\mathcal{F}$, hence their pairwise intersections have size divisible by $4$. Thus we have
\begin{align*}
&2|\supp(\bm{u})\cap\supp(\bm{v})\cap\supp(\bm{w})|\\
=&|\supp(\bm{u})\cap\supp(\bm{v})|+|\supp(\bm{v})\cap\supp(\bm{w})|-|\supp(\bm{u}+\bm{w})\cap\supp(\bm{v})|\equiv &0\pmod{4}.
\end{align*}
Note that $\supp(\bm{u}\circ \bm{v})=\supp(\bm{u})\cap \supp(\bm{v})$. The above congruent implies that $|\supp(\bm{u}\circ \bm{v})\cap \supp(\bm{w})|$ is even. Thus $(\bm{u}\circ \bm{v})\cdot \bm{w}=0$ in $\F_2$ and so $\bm{u}\circ \bm{v}\in \mc{C}^\perp$. 

Moreover, since $\wt(\bm{u}\circ\bm{v})=|\supp(\bm{u})\cap\supp(\bm{v})|\equiv0\pmod{4}$, by Lemma \ref{lem_structure_C}, we have $\bm{u}\circ\bm{v}\in\mathcal{C}$.
\end{proof}

\begin{theorem}\label{thm_m4n00_odd}
Let $n\geq 3$ be odd. Then $m_{4,n}(0,0)<2^{\frac{n-1}{2}}$.
\end{theorem}
\begin{proof}
Suppose to the contrary that $m_{4,n}(0,0)\geq 2^{\frac{n-1}{2}}$. Take the above notation for $d,m$, $\mc{F}$ and $\mc{C}$. We consider the characteristic vectors over $\F_2$. 
A non-zero vector $\bm{u}\in\mathcal{C}$ is called \emph{minimal} if there is no other non-zero $\bm{v}\in\mathcal{C}$ 
with $\supp(\bm{v})\subsetneq\supp(\bm{u})$. 
Let $\bm{u}_1,\ldots,\bm{u}_r$ be all the minimal vectors of $\mathcal{C}$.

We claim that the supports of $\bm{u}_1,\ldots,\bm{u}_r$ are pairwise disjoint. 
Indeed, for $i\neq j$, Lemma \ref{lem_Hadamard_closed} gives $\bm{u}_i\circ\bm{u}_j\in\mc{C}$, 
and $\supp(\bm{u}_i\circ\bm{u}_j)=\supp(\bm{u}_i)\cap\supp(\bm{u}_j)$. 
Minimality forces this intersection to be empty.

Hence $\bm{u}_1,\ldots,\bm{u}_r$ are linearly independent. 
Moreover, every vector in $\mathcal{C}$ is a sum of some minimal vectors. 
To see this, take any non-zero $\bm{w}\in\mc{C}$. 
The set $\{\bm{v}\in\mc{C}\setminus\{\bm{0}\}:\supp(\bm{v})\subseteq\supp(\bm{w})\}$ is non-empty 
(it contains $\bm{w}$), so we can pick a minimal vector $\bm{u}_i$ from it. 
Then $\bm{w}+\bm{u}_i\in\mc{C}$ and $\wt(\bm{w}+\bm{u}_i)<\wt(\bm{w})$. 
By induction on weight, $\bm{w}+\bm{u}_i$ is a sum of minimal vectors, and the claim follows.

Thus $\bm{u}_1,\ldots,\bm{u}_r$ form a basis for $\mc{C}$, and consequently $r=d$.
Since each minimal vector is non-zero, its weight is at least $4$. Then we have
\[
n\geq \sum_{i=1}^{d} |\supp(\bm{u}_i)| \geq 4d = 2n-2,
\]
which is impossible for $n\geq 3$. This contradiction completes the proof.
\end{proof}

Note that if $k\equiv 0\pmod{4}$, then any $(0,0)$-town$\pmod{k}$ family is also a $(0,0)$-town$\pmod{4}$ family. From the above Theorems \ref{thm_m4n00_even} and \ref{thm_m4n00_odd}, we readily obtain the following result.

\begin{corollary}\label{cor_m4n00}
Let $k\equiv 0\pmod{4}$. Then $m_{k,n}(0,0)<2^{\lfloor \frac{n}{2}\rfloor}$.
\end{corollary}

Having established the above preliminary results, we now turn to the proof of Theorem \ref{thm_main_aa}.

\begin{proof}[Proof of Theorem \ref{thm_main_aa}]
The upper bound $m_{k,n}(a,a)\leq 2^{\lfloor \frac{n}{2}\rfloor}$ has been proved in Corollary \ref{cor_maa_ub}. The sufficient condition for the equality to hold is obvious. So it remains to derive the necessity.

Suppose that $m_{k,n}(a,a)=2^{\lfloor \frac{n}{2}\rfloor}$. Then by Lemma \ref{lem:aa-le-00}, we have $m_{k,n}(0,0)=2^{\lfloor \frac{n}{2}\rfloor}$. It then follows from Proposition \ref{pp_00towm_1} and Corollary \ref{cor_m4n00} that $k=2$ or $n=k=3$. When $n=k=3$, it is easy to check that only the $(0,0)$-town$\pmod{3}$ family $\{\varnothing,\{1,2,3\}\}$ can attain the upper bound $2^{\lfloor \frac{n}{2}\rfloor}=2$. As we have introduced in Section \ref{sec_intro}, when $k=2$, we must have $a=0$, or $a=1$ and $n$ is odd. This completes the proof.
\end{proof}

\subsection{An improved upper bound for $m_{3,n}(a,a)$}\label{sub_m3naa}
In this subsection, we work over the finite field $\F_3$. Let $\mc{F}=\{F_1,\ldots,F_m\}$ be a $(0,0)$-town$\pmod{3}$ family on $[n]$, and let $\bm{v}_1,\ldots,\bm{v}_m$ be their respective characteristic vectors in $\F_3^n$. Let $V=\<\bm{v}_1,\ldots,\bm{v}_m\>$. Then $V$ is a totally isotropic subspace of $\F_3^n$.

%

\begin{theorem}\label{thm_m3n00}
Let $n\equiv 0\pmod{4}$. Then $m_{3,n}(0,0)< (\frac{32}{9})^{\frac{n}{4}}\approx 1.3732^n$.
\end{theorem}
\begin{proof}
Let $M$ be a maximal totally isotropic subspace of $\F_3^n$ containing $V$. By Remark \ref{rmk_tispace}, we have $\dim(M)=\frac{n}{2}$ and therefore $M=M^\perp$. Given $\bm{x}\in M$, we have $\wt(\bm{x})\equiv x_1^2+\cdots+x_n^2\equiv \bm{x}\cdot\bm{x}\equiv 0\pmod{3}$. Thus $M$ is a Type III code.

We will bound $m_{3,n}(0,0)$ by calculating the number of 0-1 vectors in $M$. Let $S=\{0,1\}^n$ and let $\mathds{1}_M(\bm{x})$ be the indicator function of $M$. Using \eqref{eq_indi} and substituting $M^\perp$ by $M$, we deduce that 
\begin{align}
|M\cap S|&=\sum_{\bm{x}\in S}\mathds{1}_M(\bm{x})=\frac{1}{|M|}\cdot\left\vert\sum_{\bm{u}\in M}\sum_{\bm{x}\in S}\zeta_3^{\bm{x}\cdot \bm{u}}\right\vert\nonumber\\
&=\frac{1}{|M|}\cdot\left\vert\sum_{\bm{u}\in M}\sum_{x_1,\ldots,x_n\in \{0,1\}}\zeta_3^{x_1u_1}\cdots\zeta_3^{x_nu_n}\right\vert\nonumber\\
&\leq\frac{1}{|M|}\sum_{\bm{u}\in M}\left\vert\prod_{j=1}^n(1+\zeta_3^{u_j})\right\vert.\label{eq_m0300_1}
\end{align}
Note that 
\[|1+\zeta_3^{u_j}|=\begin{cases}
2,&\text{if }u_j=0,\\
|1+\zeta_3|=|e^{\frac{\pi i}{3}}|=1,&\text{if }u_j=1,\\
|1+\zeta_3^2|=|e^{-\frac{\pi i}{3}}|=1,&\text{if }u_j=2.
\end{cases}\]
Thus, we have
\[\left\vert\prod_{j=1}^n(1+\zeta_3^{u_j})\right\vert=2^{n-\wt(\bm{u})}1^{\wt{(\bm{u}})}.\]
Continuing from \eqref{eq_m0300_1}, we deduce that
\[|M\cap S|\leq \frac{1}{|M|}\sum_{\bm{u}\in M}2^{n-\wt(\bm{u})}1^{\wt{(\bm{u}})}=\frac{1}{|M|}W_{M}(1,2)<\frac{1}{3^{\frac{n}{2}}}34^{\frac{n}{4}}=\left(\frac{34}{9}\right)^{\frac{n}{4}}.\]
The last inequality follows from Corollary \ref{cor_wt12}. This completes the proof.
\end{proof}

\begin{corollary}\label{cor_m3naa}
For $n\geq 3$, we have $m_{3,n}(a,a)<2.71(\frac{34}{9})^{\frac{n}{4}}\approx 2.71\cdot1.3732^n$ for $a\in \{0,1,2\}$.
\end{corollary}
\begin{proof}
If $n\equiv 0\pmod{4}$, the result follows directly from Lemma \ref{lem:aa-le-00} and Theorem \ref{lem:aa-le-00}. If $n\equiv r\pmod{4}$ where $r\in \{1,2,3\}$, we deduce from Lemma \ref{lem:aa-le-00} and Theorem \ref{lem:aa-le-00} that $m_{3,n}(0,0)\leq m_{3,n+4-r}(0,0)<(\frac{34}{9})^{\frac{n+4-r}{4}}<2.71(\frac{34}{9})^{\frac{n}{4}}$.
\end{proof}

Recall that Theorem \ref{thm_main_aa} gives the general bound $m_{3,n}(a,a)<2^{\lfloor n/2\rfloor}\approx 1.4142^n$. Corollary \ref{cor_m3naa} improves the exponential base to  $1.3732$.

\section{The case $a\not\equiv b\pmod{k}$}\label{sec_ab}

\subsection{The general upper bound for $m_{k,n}(a,b)$}

This subsection is devoted to the proof of Theorem \ref{thm_main_ab}, i.e., $m_{k,n}(a,b)\leq n$ whenever $a\not\equiv b\pmod{k}$. Throughout this subsection, we always assume that $a,b,k,n$ are integers such that $n\geq k\geq 2$, $0\leq a,b\leq k-1$ and $a\not\equiv b\pmod{k}$. Necessarily, there must be some prime $p$ and a positive integer $t$ such that $p^t\mid k$ and $p^t\nmid a-b$. 

%
%
%
%
Our proof consists of two steps. We will first show that an $(a,b)$-town$\pmod{k}$ family has size at most $n+1$, and then we show that the upper bound $n+1$ can never be achieved. 

\begin{theorem}\label{thm_ab_n+1}
Take the above notation. Let $\mc{F}=\{F_1,\ldots,F_m\}$ be an $(a, b)$-town$\pmod{k}$ family on $[n]$, where $a \not\equiv b \pmod{k}$. Then $|\mathcal{F}| \leq n+1$.
\end{theorem}
\begin{proof} 
Take $f\in\Z$ such that $b+kf<0$ and $-(b+kf)$ is not a perfect square. By Lagrange's Four Square Theorem, there exists $z_1,\ldots,z_N\in\Z^+$ with $2\leq N\leq 4$ such that 
\[-b-kf=z_1^2+\cdots+z_N^2.\]

For each $F_i\in\mc{F}$, define its extended characteristic vector in $\mathbb{Q}^{n+N}$ by
\[(\bm{\chi}_i)_j=\begin{cases}
1,&\text{if $1\leq j\leq n$ and $j\in F_i$},\\
0,&\text{if $1\leq j\leq n$ and $j\notin F_i$},\\
z_s,&\text{if }j=n+s, 1\leq s\leq N.
\end{cases}\]
That is, $\bm{\chi}_i=(\bm{v}_i,z_1,\ldots,z_N)$, where $\bm{v}_i$ is the characteristic vector of $F_i$. Then we have
\begin{equation}\label{eq_chichi}
\bm{\chi}_i\cdot \bm{\chi}_j=|F_i\cap F_j|+z_1^2+\cdots+z_N^2\equiv
\begin{cases}
a-b\pmod{k},&\text{if }i=j,\\
0\pmod{k},&\text{if }i\neq j.
\end{cases}
\end{equation}
It follows that 
\begin{equation}\label{eq_chichi2}
\bm{\chi}_i\cdot \bm{\chi}_j\equiv
\begin{cases}
a-b\not\equiv 0\pmod{p^t},&\text{if }i=j,\\
0\pmod{p^t},&\text{if }i\neq j.
\end{cases}
\end{equation}
Suppose that $c_1\bm{\chi}_1+\cdots+c_m\bm{\chi}_m=0$ for some $c_1,\ldots,c_m\in \mathbb{Q}$. If not all the $c_i$'s are zero, by clearing denominators and then dividing by the greatest common divisor of all resulting integer coefficients, we obtain a primitive relation
\begin{equation}\label{eq_chilb_0}
c_1'\bm{\chi}_1+\cdots+c_m'\bm{\chi}_m=0,
\end{equation}
where $c_i'\in\Z$ and $\gcd(c_1',\ldots,c_m')=1$. 

For any $i\in[m]$, we deduce from \eqref{eq_chilb_0} and \eqref{eq_chichi2} that
\[c_i'(\bm{\chi}_i\cdot \bm{\chi}_i)\equiv 0\pmod{p^t}.\]
Since $\bm{\chi}_i\cdot\bm{\chi}_i\not\equiv 0\pmod{p^t}$, we conclude that $p\mid c_i'$. This contradicts the fact that $\gcd(c_1',\ldots,c_m')=1$. Therefore $c_1=\cdots=c_m=0$, and consequently the vectors $\bm{\chi}_1,\ldots,\bm{\chi}_m$ are linearly independent. 

Define $N-1$ vectors in $\Q^{n+N}$ as follows.
\begin{align*}
\bm{\omega}_2=&(0,\ldots,0,z_2,-z_1,0,\ldots,0),\\
\bm{\omega}_3=&(0,\ldots,0,z_3,0,-z_1,0,\ldots,0),\\
\vdots\\
\bm{\omega}_N=&(0,\ldots,0,z_N,0,\ldots,0,-z_1).
\end{align*}
That is, for each $i\in\{2,\ldots,N\}$, the $(n+1)$-th position of $\bm{\omega}_{i}$ is $z_i$, the $(n+i)$-th position of $\bm{\omega}_i$ is $-z_1$, and the rest positions are 0. 
Then $\bm{\chi}_1,\ldots,\bm{\chi}_m$ belong to the subspace $\<\bm{\omega}_2,\ldots,\bm{\omega}_N\>^\perp$ of $\mathbb{Q}^{n+N}$. Since $\dim(\<\bm{\omega}_2,\ldots,\bm{\omega}_N\>)=N-1$, we conclude that 
\[m\leq \dim(\<\bm{\omega}_2,\ldots,\bm{\omega}_N\>^\perp)=n+N-(N-1)=n+1.\] 
This completes the proof.
\end{proof}

Next, we show that the upper bound $n+1$ can not be reached.

\begin{proof}[Proof of Theorem \ref{thm_main_ab}]
Take the notation from Theorem \ref{thm_ab_n+1} and its proof. If $\<\bm{v}_{1},\ldots,\bm{v}_m\>\neq \mathbb{Q}^n$, then we can choose a non-zero vector $\bm{\eta}\in\<\bm{v}_{1},\ldots,\bm{v}_m\>^\perp$. Then we define
\[\widetilde{\bm{\eta}}=(\bm{\eta},0,\ldots,0)\in \Q^{n+N}.\]
Thus $\bm{\chi}_1,\ldots,\bm{\chi}_m$ belong to the subspace $\<\widetilde{\bm{\eta}},\bm{\omega}_2,\ldots,\bm{\omega}_N\>^\perp$ of $\Q^{n+N}$. Since $\dim(\<\widetilde{\bm{\eta}},\bm{\omega}_2,\ldots,\bm{\omega}_N\>)=N$, we have $m\leq (n+N)-N=n$. 

Next, we deal with the case when $\<\bm{v}_{1},\ldots,\bm{v}_m\>= \mathbb{Q}^n$. Recall from Theorem \ref{thm_ab_n+1} that $m\leq n+1$. Suppose to the contrary that $m=n+1$ and assume with loss of generality that $\bm{v}_{1},\ldots,\bm{v}_{n}$ form a basis for $\Q^n$. Then there exist integers $\mu_0,\mu_1,\ldots,\mu_n$ with $\gcd(\mu_0,\mu_1, \dots, \mu_n) = 1$ such that
\begin{equation}\label{eq_bm1_0}
\mu_0\bm{1} = \mu_1 \bm{v}_{1} + \cdots + \mu_n \bm{v}_{n}.
\end{equation}
Computing the inner product of both sides of \eqref{eq_bm1_0} with $\bm{v}_{n+1}$, we have
\begin{equation}\label{eq_mu_1}
\mu_0a\equiv \sum_{j=1}^n\mu_jb\pmod{k}.
\end{equation}

For each $i\in [n]$, computing the inner product of both sides of \eqref{eq_bm1_0} with $\bm{v}_i$, we have
\begin{equation}\label{eq_mu_2}
\mu_0 a\equiv \mu_i a+\sum_{\substack{j=1\\ j\neq i}}^n \mu_j b \pmod{k}.
\end{equation}
Combining \eqref{eq_mu_1} and \eqref{eq_mu_2}, we obtain
\[\mu_i(a-b)\equiv 0\pmod{k}.\]
Since $p^t\mid k$ and $p^t\nmid a-b$, we deduce that
\[p\mid \mu_i, \forall i\in[n].\]
Since $\bm{v}_1,\ldots,\bm{v}_n$ are 0-1 vectors, it then follows from \eqref{eq_bm1_0} that $p\mid \mu_0$. But this contradicts the assumption that $\gcd(\mu_0,\mu_1,\ldots,\mu_n)=1$. 
This completes the proof.
\end{proof}

Theorem \ref{thm_main_ab} gives the upper bound $m_{k,n}(a,b)\leq n$ for all $a\not\equiv b\pmod{k}$.  The following corollary exhibits infinitely many cases where this bound is attained.

\begin{corollary}\label{cor_pg2q}
Let $q$ be a prime power and let $k\geq 2$ be such that $q\not\equiv 0\pmod{k}$. Then $m_{k,q^2+q+1}(q+1,1)=q^2+q+1$.
\end{corollary}
\begin{proof}
Let $\mc{P}$ and $\mc{L}$ be the point set and the line set of the projective plane $\PG(2,q)$, respectively. We have $|\mc{P}|=|\mc{L}|=q^2+q+1$. Note that each line contains $q+1$ points and any two lines intersect in exactly one point. Then $\mc{L}$ forms a $(q+1,1)$-town$\pmod{k}$ family on $\mc{P}$ for any $k\geq 2$, giving the lower bound $m_{k,q^2+q+1}(q+1,1)\geq q^2+q+1$.

If $q\not\equiv 0\pmod{k}$, then $q+1\not\equiv 1\pmod{k}$, so Theorem \ref{thm_main_ab} applies and yields the matching upper bound. Hence equality holds. 
\end{proof}

\subsection{Exact values and improved bounds in special cases}\label{sub_betterab}

\subsubsection{The maximum size of $(t, t-1)$-town$\pmod{k}$ families for general $k$}

\begin{theorem}\label{thm_mkntt-1}
Let $n\equiv t+1\pmod{k}$ where $1\leq t\leq k$. Then $m_{k,n}(t,t-1)=n$.
\end{theorem}
\begin{proof}
By Proposition \ref{prop_basic}(iii), we have $m_{k,n}(t,t-1)=m_{k,n}(n-t,n-2t+t-1)$. Since $n\equiv t+1\pmod{k}$, then we have $m_{k,n}(n-t,n-2t+t-1)=m_{k,n}(1,0)=n$. The last equality holds because $\{\{1\},\ldots,\{n\}\}$ is a $(1,0)$-town$\pmod{k}$ family of size $n$. 
\end{proof}

\begin{remark}
When $t=1$, $m_{k,n}(1,0)=n$ holds for all $n$ since the upper bound is trivially attained by the family of all singletons.
\end{remark}

\begin{theorem}\label{thm_tt1n1}
Let $n\equiv t\pmod{k}$ where $2\leq t\leq k$. Then $m_{k,n}(t,t-1)\leq n-1$.
\end{theorem}
\begin{proof}
Suppose to the contrary that $m_{k,n}(t,t-1)=n$ and suppose that $\{F_1,\ldots,F_n\}$ is a $(t,t-1)$-town$\pmod{k}$ family. Let $\bm{v}_1,\ldots,\bm{v}_n$ be their respective characteristic vectors. Plugging $(a,b,k,m)=(t,t-1,k,n)$ in Lemma \ref{lem_linindep}, we have
\[gcd(a(a-b),(a-b)(a-b-mb),k)=\gcd(t,t^2-t+1,k)=1<k,\]
so $\bm{v}_1,\ldots,\bm{v}_n$ are linearly independent over $\Q$. Then there exist $b_1,\ldots,b_n\in \Q$ such that $b_1\bm{v}_1+\cdots+b_n\bm{v}_n=\bm{1}$. By clearing denominators and then dividing by the greatest common divisor of all resulting integer coefficients, we obtain a primitive relation
\begin{equation}\label{eq_tt1ee1}
\mu_1\bm{v}_1+\cdots+\mu_n\bm{v}_n=\mu_0\bm{1},
\end{equation}
where $\mu_0,\ldots,\mu_n\in\Z$ and $\gcd(\mu_0,\ldots,\mu_n)=1$. 

For each $i\in[n]$, taking inner products of both sides of \eqref{eq_tt1ee1} with $\bm{v}_i$, we have
\begin{equation}\label{eq_tt1ee2}
t\mu_i+(t-1)\sum_{\substack{j=1\\j\neq i}}^n\mu_j\equiv t\mu_0.
\end{equation}
This implies that
\[\mu_i\equiv t\mu_0-(t-1)\sum_{j=1}^n\mu_j\equiv S\pmod{k},\]
for some $S\in\Z$. Thus \eqref{eq_tt1ee2} is reduced to
\begin{equation}\label{eq_tt1ee3}
tS+(t-1)(n-1)S\equiv (t^2-t+1)S\equiv t\mu_0\pmod{k}.
\end{equation}

Taking inner products of both sides of \eqref{eq_tt1ee1} with $\bm{1}$, we have
\begin{equation*}
t\sum_{j=1}^n\mu_j\equiv n\mu_0\equiv t\mu_0\pmod{k}.
\end{equation*}
Together with \eqref{eq_tt1ee3}, we obtain
\begin{equation}\label{eq_tt1ee4}\
(t-1)S\equiv 0\pmod{k}.
\end{equation}
Let $d=\gcd(t-1,k)$. Then $d<k$, and \eqref{eq_tt1ee4} implies that $\frac{k}{d}\mid \mu_i$ for each $i\in [n]$. By \eqref{eq_tt1ee1}, we see that $\frac{k}{d}$ divides each position of $\mu_1\bm{v}_1+\cdots+\mu_n\bm{v}_n$ and so $\frac{k}{d}$ divides each position of $\mu_0\bm{1}$, that is, $\frac{k}{d}\mid \mu_0$. This contradicts the assumption that $\gcd(\mu_0,\ldots,\mu_n)=1$. Therefore $m_{k,n}(t,t-1)\leq n-1$. 
\end{proof}

Together with Example \ref{eg_ab}, we can determine $m_{k,n}(2,1)$ when $n\equiv 2\pmod{k}$.
\begin{corollary}\label{cor_mkn21}
Let $n\equiv 2\pmod{k}$. Then $m_{k,n}(2,1)=n-1$.
\end{corollary}

\subsubsection{The maximum size of $(0, p-1)$-town$\pmod{p}$ families for odd prime $p$}

\begin{theorem}\label{thm_mpn0p-1}
Let $p$ be an odd prime and $n\equiv 0\pmod{p}$. Then $m_{p,n}(0,p-1)\leq n-2$.
\end{theorem}
\begin{proof}
By Theorem \ref{thm_tt1n1}, we have $m_{p,n}(0,p-1)\leq n-1$. Suppose that $m_{p,n}(0,p-1)=n-1$ and assume that $F_1,\ldots,F_{n-1}$ form a $(0,p-1)$-town$\pmod{p}$ family. Let $\bm{v}_1,\ldots,\bm{v}_{n-1}$ be their respective characteristic vectors in $\F_p$. Then by Lemma \ref{lem_linindep} with $(a,b,k,m)=(0,p-1,p,n-1)$, we have
\[\gcd(a(a-b),(a-b)(a-b+mb),k)=\gcd(2,p)<p,\]
and so $\bm{v}_1,\ldots,\bm{v}_{n-1}$ are linearly independent over $\Q$.

Note that $\bm{1}\cdot \bm{v}_i\equiv 0\pmod{p}$ for each $i\in [n-1]$. Then $\bm{v}_1,\ldots,\bm{v}_{n-1}$ form a basis for $\bm{1}^\perp$. Since $n\equiv 0\pmod{p}$, then $\bm{1}\in \bm{1}^\perp=\<\bm{v}_1,\ldots,\bm{v}_{n-1}\>$. So we have
\begin{equation}\label{eq_mpn0p-1_1}
\mu_1\bm{v}_1+\cdots+\mu_{n-1}\bm{v}_{n-1}=\mu_0\bm{1},
\end{equation}
where $\mu_0,\mu_1,\ldots,\mu_{n-1}\in\Q$ and $\gcd(\mu_0,\ldots,\mu_{n-1})=1$. For each $i\in[n-1]$, we have 
\[(p-1)\sum_{\substack{j=1\\j\neq i}}^{n-1}\mu_j\equiv 0\pmod{p},\]
which implies that $\mu_j\equiv \sum_{j=1}^{n-1}\mu_j\equiv S\pmod{p}$ for some $S\in \Z$, and so 
\[0\equiv \sum_{\substack{j=1\\j\neq i}}^{n-1}\mu_j\equiv (n-2)S\equiv -2S\pmod{p}.\]
Since $p$ is odd, then $S\equiv 0\pmod{p}$. By \eqref{eq_mpn0p-1_1}, this implies that $\mu_0\equiv 0\pmod{p}$, a contradiction. Therefore $m_{p,n}(0,p-1)\leq n-2$. 
\end{proof}

\begin{corollary}\label{cor_m3n02}
Let $n\equiv 0\pmod{3}$. Then $m_{3,n}(0,2)=n-2$.
\end{corollary}
\begin{proof}
Letting $p=3$ in Theorem \ref{thm_mpn0p-1} yields $m_{3,n}(0,2)\leq n-2$. Example \ref{eg_ab} implies that there is a $(0,2)$-town$\pmod{3}$ family on $[n]$ of size $n-2$. Therefore $m_{3,n}(0,2)=n-2$.
\end{proof}

\subsubsection{The maximum size of $(t,t-1)$-town$\pmod{4}$ families}
Let $n=4c$ for some $c\in \Z^+$, so $n\equiv 0\pmod{4}$. Let $\{F_1,\ldots,F_m\}$ be a $(0,3)$-town$\pmod{4}$ family on $[n]$ and let $\bm{v}_1,\ldots,\bm{v}_m$ be their respective characteristic vectors in $\F_2^n$.

\begin{lemma}\label{lem_m4n03}
If $m=n-2$, then the vectors $\bm{v}_1,\ldots,\bm{v}_m$ are linearly independent over $\F_2$, and $\bm{1}\notin\<\bm{v}_1,\ldots,\bm{v}_m\>_{\F_2}$.
\end{lemma}
\begin{proof}
It is clear that 
\[\bm{v}_i\cdot \bm{v}_j\equiv \begin{cases}
0\pmod{2}, &\text{if }i=j,\\
1\pmod{2},&\text{if }i\neq j.\end{cases}\]
View $\bm{v}_1,\ldots,\bm{v}_m$ as column vectors and define a matrix $A=(\bm{v}_1\mid \ldots\mid\bm{v}_m)$. Then 
\[B:=A^\top A\equiv J+I\pmod{2}.\]
Note that
\[\det(B)=m+1\equiv 1\pmod{2}.\]
Thus $B$ is invertible over $\F_2$ and so $\rank_{\F_2}(A)=m$. This implies that $\bm{v}_1,\ldots,\bm{v}_m$ are linearly independent over $\F_2$.

To prove the second statement, we suppose to the contrary that $\mu_1\bm{v}_1+\cdots+\mu_m\bm{v}_m=\bm{1}$ for some $\mu_1,\ldots,\mu_n\in\F_2$. For $i\in[m]$, taking inner product of both sides of the above equation with $\bm{v}_i$, we obtain $\sum_{j\neq i}\mu_j=0$. It follows that $\mu_i=\sum_{j=1}^m\mu_j=S$ for some $S\in\F_2$. Thus $0=\sum_{j\neq i}\mu_j=(m-1)S=S$. It follows that $\bm{0}=\bm{1}$, a contradiction. Therefore $\bm{1}\notin\<\bm{v}_1,\ldots,\bm{v}_m\>_{\F_2}$.
\end{proof}

Let $V=\F_2^n$, and let $W$ be the subspace of $V$ consisting of vectors of even weight. Define a quadratic form on $W$ by 
\[Q(\bm{x})\equiv \frac{\wt(\bm{x})}{2}\pmod{2}.\]

Let $U=\<\bm{v}_1,\ldots,\bm{v}_m\>_{\F_2}$. We count the number of singular vectors in $W$ and $U$, respectively. Define
\begin{align*}
N_W&=|\{\bm{x}\in W: Q(\bm{x})=0\}|,\\
N_U&=|\{\bm{x}\in U: Q(\bm{x})=0\}|.
\end{align*}

\begin{lemma}\label{lem_NWNU}
If $m=n-2$, then $N_W=2^{4c-2}+2^{2c-1}(-1)^{c}$ and $N_U=2^{4c-3}+2^{2c-2}(-1)^{c-1}$. 
\end{lemma}
\begin{proof}
We first compute $N_W$. For $\bm{x}\in W$, it is clear that $Q(\bm{x})=0$ if and only if $\wt(\bm{x})\equiv 0\pmod{4}$. Thus
\[N_W=\sum_{\substack{j=0\\ j\equiv 0\pmod{4}}}^n \binom{n}{j}.\]
By Lemma \ref{lem_Tmr4}, we have
\begin{align*}
N_W&=T(n,0,4)=2^{n-2}+2^{\frac{n-2}{2}}\cos\left(\frac{n\pi}{4}\right)=2^{4c-2}+2^{2c-1}(-1)^{c}.
\end{align*}

Next, we count the number $N_U$ of vectors $\bm{x}\in U$ such that $Q(\bm{x})=0$. For any $\bm{x}\in U$, it can be written uniquely as $\bm{x}=\sum_{i\in S}\bm{v}_i$ where $S\subseteq [m]$. Thus
\[Q(\bm{x})\equiv \sum_{i\in S}Q(\bm{v}_i)+\sum_{\substack{i,j\in S\\ i<j}}\bm{v}_i\cdot \bm{v}_j\equiv \binom{|S|}{2}\pmod{2}.\]
Note that $\binom{|S|}{2}\equiv 0\pmod{2}$ if and only if $|S|\equiv 0,1\pmod{4}$. It follows that
\[N_U=\sum_{j\equiv 0,1\pmod{4}}\binom{m}{j}.\]
Recall that $m=n-2=4c-2$. By Lemma \ref{lem_Tmr4}, we obtain
\begin{align*}
N_U&=T(m,0,4)+T(m,1,4)\\
&=2^{m-2}+2^{\frac{m-2}{2}}\cos\left(\frac{(4c-2)\pi}{4}\right)+2^{m-2}+2^{\frac{m-2}{2}}\cos\left(\frac{(4c-4)\pi}{4}\right)\\
&=2^{4c-3}+2^{2c-2}(-1)^{c-1}.
\end{align*}
This completes the proof.
\end{proof}

\begin{theorem}\label{thm_m4n03}
Let $n\equiv 0\pmod{4}$. Then $m_{4,n}(0,3)=n-3$.
\end{theorem}
\begin{proof}
Assume that $m=n-2$ and take the notation in the above two lemmas and their proofs. Note that
$W=\{\bm{x}\in V: x_1+\cdots+x_n=0\}=\bm{1}^\perp$ and $\dim(W)=n-1$. Since $\bm{v}_i\cdot \bm{1}=0$ and $\bm{1}\cdot \bm{1}=0$, we have $U\subseteq \bm{1}^\perp$ and $\bm{1}\subseteq\bm{1}^\perp$. By Lemma \ref{lem_m4n03}, we have $\dim(U)=m=n-2$ and $\bm{1}\notin U$. Thus
\[W=U\oplus \<\bm{1}\>_{\F_2}.\]
Note that for $\bm{x}\in W$, we have $Q(\bm{x}+\bm{1})=Q(\bm{x})$. This implies that there is a bijection between singular vectors in $U$ and singular vectors in $W\setminus U$. Thus $N_W=2N_U$. However, by Lemma \ref{lem_NWNU}, we have $2N_U=2^{4c-2}+2^{2c-1}(-1)^{c-1}\neq N_W$, a contradiction. Thus $m_{4,n}(0,3)<n-2$.

On the other hand, by Example \ref{eg_ab}, we see that a $(0,3)$-town$\pmod{4}$ family of size $n-3$ exists. Therefore $m_{4,n}(0,3)=n-3$. 
\end{proof}

\begin{remark}\label{rmk_m4n01}
By the same arguments above, one can prove that $m_{4,n}(0,1)\leq n-3$ if $n\equiv 0\pmod{4}$. However, the equality does not always hold. For example, by Magma \cite{Magma}, we have $m_{4,4}(0,1)=1$, $m_{4,8}(0,1)=2$, and $m_{4,12}(0,1)=9$. 
\end{remark}

\begin{corollary}\label{cor_m4n32}
Let $n\equiv 3\pmod{4}$. Then $m_{4,n}(3,2)=n-2$.
\end{corollary}
\begin{proof}
Let $\{F_1,\ldots,F_m\}$ be a $(3,2)$-town$\pmod{4}$ family on $[n]$. For $i\in [m]$, define
\[F_i'=F_i\cup \{n+1\}.\]
Then $F_1',\ldots,F_m'$ form a $(0,3)$-town$\pmod{4}$ family on $[n+1]$. Since $n+1\equiv 0\pmod{4}$, by Theorem \ref{thm_m4n03}, we have $m<(n+1)-2=n-1$.
On the other hand, Example \ref{eg_ab} implies that $m\geq n-2$. Therefore $m_{4,n}(3,2)=n-2$.
\end{proof}

\subsection{Summary of exact values and improved bounds}

We conclude this section by collecting the exact values and improved upper bounds for $m_{k,n}(a,b)$ obtained throughout Sections \ref{sec_aa} and \ref{sec_ab}. For the reader's convenience, we list them in Table \ref{tab_results}.

\begin{table}[H]
\centering
\caption{Improved bounds and exact values for specific parameters}
\label{tab_results}
\begin{tabular}{c|c|c}
\toprule
$m_{k,n}(a,b)$ & Condition & Reference \\
\midrule
$m_{3,n}(a,a)<2.71\cdot 1.3732^n$ & $k=3$ & Corollary \ref{cor_m3naa} \\
$m_{k,n}(t,t-1)=n$ & $n\equiv t+1\pmod{k}$, $1\leq t\leq k$ & Theorem \ref{thm_mkntt-1} \\
$m_{k,n}(t,t-1)\leq n-1$ & $n\equiv t\pmod{k}$, $2\leq t\leq k$ & Theorem \ref{thm_tt1n1} \\
$m_{k,n}(2,1)=n-1$ & $n\equiv 2\pmod{k}$ & Corollary \ref{cor_mkn21} \\
$m_{3,n}(0,2)=n-2$ & $n\equiv 0\pmod{3}$ & Corollary \ref{cor_m3n02} \\
$m_{4,n}(3,0)\leq n-1$ & all $n$ & Lemma \ref{lem_m4n30n} \\
$m_{4,n}(0,3)=n-3$ & $n\equiv 0\pmod{4}$ & Theorem \ref{thm_m4n03} \\
$m_{4,n}(0,1)\leq n-3$ & $n\equiv 0\pmod{4}$ & Remark \ref{rmk_m4n01}\\
$m_{4,n}(3,2)=n-2$ & $n\equiv 3\pmod{4}$ & Corollary \ref{cor_m4n32} \\
\bottomrule
\end{tabular}
\end{table}

\section{Further discussion}\label{sec_discuss}

\subsection{The parity of $m_{k,n}(a,a)$}

While the exact value of $m_{k,n}(a,a)$ seems very difficult to determine in general, a more accessible question is its parity. Based on computational evidence for small parameters, we would like to propose the following conjecture.

\begin{conjecture}
If $n\geq k+a$, then $m_{k,n}(a,a)$ is even.
\end{conjecture}

We are able to prove the following special case of this conjecture.

\begin{proposition}
If $n\geq k$ and $n\equiv 0\pmod{k}$, then $m_{k,n}(0,0)$ is even.
\end{proposition}
\begin{proof}
Let $\mc{F}$ be a $(0,0)$-town$\pmod{k}$ family with $|\mc{F}|=m_{k,n}(0,0)$. For any $F\in \mc{F}$, we have $|\overline{F}|=n-|F|\equiv 0\pmod{k}$ and $|\overline{F}\cap F'|=|F'|-|F\cap F'|\equiv 0\pmod{k}$ for any $F'\in \mc{F}$. Since $\mc{F}$ has the maximum size, we see that $\overline{F}\in\mc{F}$. Therefore $|\mc{F}|$ must be even.
\end{proof}

\subsection{Better upper bounds for $m_{k,n}(a,a)$}
For $k=3$, Theorem \ref{cor_m3naa} gives the upper bound 
$m_{3,n}(a,a)<2.71\cdot1.3732^n$. On the other hand, using 
Magma \cite{Magma} we can confirm that $m_{3,12}(0,0)=24$. 
Lemma \ref{lem_combine} then yields the lower bound 
$m_{3,n}(0,0)\geq 24^{\lfloor n/12\rfloor}\approx 1.303^n$, 
which coincides with the bound \eqref{eq_FObound} obtained by Frankl and Odlyzko 
via Hadamard matrices. It remains an interesting open question to determine the true
exponential growth rate of $m_{3,n}(0,0)$.

By Proposition~\ref{prop_basic}(iii), we have $m_{k,n}(k-1,k-1)=m_{k,n}(0,0)$ whenever $n\equiv k-1\pmod{k}$. The case $n\not\equiv k-1\pmod{k}$ is less clear. For $k=2$, Theorem~\ref{thm_main_aa} gives $m_{2,n}(1,1)<m_{2,n}(0,0)$ for all even $n$. The following theorem covers the case $n\equiv0\pmod{k}$ for general $k$.

\begin{theorem}\label{thm:strict}
If $n\equiv 0\pmod{k}$ and $1\le a\le k-1$, then 
$m_{k,n}(a,a)<m_{k,n}(0,0)$.
\end{theorem}

\begin{proof}
Suppose that $\mathcal{F}$ is an $(a,a)$-town$\pmod{k}$ family 
with $|\mathcal{F}|=m_{k,n}(0,0)$. Fix $S\in\mathcal{F}$ and define
\[
\mathcal{G}=\{S\triangle A: A\in\mathcal{F}\}.
\]
Then $\mathcal{G}$ is a $(0,0)$-town$\pmod{k}$ family of the 
maximum size. Since $n\equiv0\pmod{k}$, we have $[n]\in\mathcal{G}$. 
It follows that $\overline{S}:=[n]\setminus S\in\mathcal{F}$, which 
leads to a contradiction because $|S\cap\overline{S}|=0$. Therefore 
$|\mathcal{F}|<m_{k,n}(0,0)$.
\end{proof}

For the remaining residue classes $n\not\equiv0,k-1\pmod{k}$, the following problem remains open.

\begin{problem}\label{prob:compare}
Let $n\not\equiv k-1\pmod{k}$. Is it true that 
$m_{k,n}(k-1,k-1)<m_{k,n}(0,0)$?
\end{problem}

\subsection{The special case \(m_{k,n}(t,t-1)\)}

When $n\equiv t\pmod{k}$ where $1\leq t\leq k$, Example \ref{eg_ab} always supplies a $(t,t-1)$-town$\pmod{k$} family of size $n-t+1$, so
\[
m_{k,n}(t,t-1)\geq n-t+1\qquad\text{for } n\equiv t\pmod{k}.
\]
The classical results for $k=2$ (see Section \ref{sec_intro}) together with the exact values
obtained in Subsection \ref{sub_betterab} show that this lower bound is tight for all $t$ when
$k=2,3,4$, and for $t=1,2$ with any modulus $k$. These results suggest that the
lower bound might always be tight when $n\equiv t\pmod{k}$.

\begin{problem}\label{prob:t-t-1}
Let \(n\equiv t\pmod{k}\) with \(1\le t\le k\).  Is it true that
\[
m_{k,n}(t,t-1)=n-t+1?
\]
Equivalently, by taking complements (Proposition \ref{prop_basic}(iii)), is it true that
\[
m_{k,n}(0,k-1)=n-t+1?
\]
\end{problem}

\subsection{The example from $\AG(2,q)$}

We have seen in Example \ref{eg_affine} that the affine plane $\AG(2,q)$ induces a $(q+1,1)$-town$\pmod{k}$ family of size $q^2$ on a $q(q+1)$-set for any $k\geq 2$. In other words, we have
\[m_{k,q(q+1)}(q+1,1)\geq q^2.\]
We check by Magma \cite{Magma} that the equality holds for $q=2, k\geq 3$, and $q=3, k\geq 4$. This motivates the following question.
\begin{problem}
Let $q$ be a prime power. Is it true that $m_{k,q(q+1)}(q+1,1)=q^2$ if $k\geq q+1$?
\end{problem}

\section*{Acknowledgement} 
The author acknowledges the support of the National Natural Science Foundation of China (grant No. 12461061), and the support of the Yunnan Fundamental Research Projects (grant No. 202601AT070176).

\bibliographystyle{plain}

\begin{thebibliography}{77}
\bibitem{AC22} M. Antipov and D. Cherkashin, Lov\'{a}sz theta approach to eventown problem, {\it Linear Algebra Appl.} {\bf 655} (2022), 302--313.

\bibitem{Ber69} E. R. Berlekamp, On subsets with intersections of even cardinality, {\it Canad. Math. Bull.} {\bf 12} (1969), 471--474.

\bibitem{BF22} L. Babai and P. Frankl, {\it Linear algebra methods in combinatorics}, Preliminary Version 2.2, Department of Computer Science, University of Chicago, 2022. 

\bibitem{Magma} W. Bosma, J. Cannon, C. Fieker, A. Steel, {\it Handbook of Magma Functions}, 2017.

\bibitem{BCZ25}B. Bukh, T.-W. Chao, Z. Zheng, The Oddtown problem modulo a composite number, arXiv:2509.00586, 2025.

\bibitem{FO83} P. Frankl and A. M. Odlyzko, On subsets with cardinalities of intersections divisible by a fixed integer, {\it European J. Combin.} {\bf 4} (1983), 215--220.

\bibitem{FT16} P. Frankl and N. Tokushige, Uniform eventown problems, {\it European J. Combin. } {\bf 51} (2016), 280--286.

\bibitem{GST22} L. Gishboliner, B. Sudakov, I. Tomon, Small doubling, atomic structure and $\ell$-divisible set families, {\it Discrete Anal.} 2022, 11.


\bibitem{Gra75} J. E. Graver, Boolean designs and self-dual matroids, {\it Linear Algebra Appl.} {\bf 10} (1975), 111--128.



\bibitem{HP03} W.C. Huffman and V. Pless, {\it Fundamentals of error-correcting codes}, Cambridge University Press, Cambridge, 2003

\bibitem{JON25} G. Johnston and J. O'Neill, A few new oddtown and eventown problems, {\it J. Comb.} {\bf 16} (2025), 303--326.

\bibitem{Lin26} N. Lindzey, An eventown result for permutations, arXiv:2601.12613, 2026. 

\bibitem{Odl81} A. M. Odlyzko, On the ranks of some (0,1)-matrices with constant row sums.
{\it J. Austral. Math. Soc. Ser. A} {\bf 31} (1981), 193--201.

\bibitem{ONe23} J. O'Neill, A short note on supersaturation for oddtown and eventown, {\it Discrete Appl. Math.} {\bf 334} (2023), 63--67.

\bibitem{ONV22} J. O'Neill and J. Verstra\"{e}te, A note on $k$-wise oddtown problems, {\it Graphs Combin.} {\bf 38} (2022), 101. 

\bibitem{Ram1834} C. Ramus, Solution g\'{e}n\'{e}rale d\textquotesingle un probl\`{e}me d\textquotesingle analyse combinatoire, {\it J. Reine Angew. Math.} {\bf 11} (1834), 353--355.

\bibitem{SV18} B. Sudakov and P. Vieira, Two remarks on eventown and oddtown problems, {\it SIAM J. Discrete Math.} {\bf 32} (2018), 280--295.


\bibitem{WZG25} X. Wei, X. Zhang, G. Ge, On set systems with strongly restricted intersections, {\it Des. Codes Cryptogr.} {\bf 93} (2025), 667--682.

\bibitem{WZZG25} X. Wei, Y. Zhao, X. Zhang, G. Ge, On supersaturation for oddtown and eventown, {\it SIAM J. Discrete Math.} {\bf 39} (2025), 2311--2327.

\bibitem{VM25} N. Veselinov and M. Marinov, On the maximal size of $(a,b)$-town$\pmod{k}$ families, arXiv:2510.00251, 2025.

\end{thebibliography}

\end{document}